\documentclass[smallextended,referee,envcountsect,]{svjour3}
\smartqed
\usepackage{graphicx}
\usepackage{mathrsfs}
\usepackage{algorithm,algorithmic}
\usepackage{amsmath,amssymb}

\usepackage{multirow}
\usepackage{subfigure}

\newtheorem{mytheo}{Theorem}
\newtheorem{mydef}{Definition}

\newtheorem{myrema}{Remark}
\newtheorem{mylemma}{Lemma}
\newtheorem{myassump}{Assumption}

\usepackage{color}
\usepackage{bm}
\usepackage{hyperref}
\journalname{JOTA}

\begin{document}

\title{Approximate Methods for Solving Chance Constrained Linear Programs in Probability Measure Space}

%\subtitle{Using  the  LaTex Template}

\author{Xun Shen         \and
        Satoshi Ito %etc.
}

\institute{ X. Shen \at
              Graduate School of Engineering, \\
              Osaka University \\
              Osaka 565-0871 Japan\\
              \email{shenxun@eei.eng.osaka-u.ac.jp}           %  \\
           \and
           S. Ito \at
              Department of Statistical Inference and Mathematics \\
              The Institute of Statistical Mathematics \\
              Tokyo 190-8562, Japan\\
              \email{sito@ism.ac.jp}
}

\date{Received: date / Accepted: date}
%The correct dates will be entered by the editor.

\maketitle

\begin{abstract}
A risk-aware decision-making problem can be formulated as a chance-constrained linear program in probability measure space. Chance-constrained linear program in probability measure space is intractable, and no numerical method exists to solve this problem. This paper presents numerical methods to solve chance-constrained linear programs in probability measure space for the first time. We propose two solvable optimization problems as approximate problems of the original problem. We prove the uniform convergence of each approximate problem. Moreover, numerical experiments have been implemented to validate the proposed methods.
\end{abstract}

%All acknowledgements should be placed in the back of the paper after Conclusions..

\section{Introduction}
\label{intro}
Let $\mathcal{X}\subset\mathbb{R}^n$ be a compact set with the infinity norm defined by $\|x\|_{\infty}=\max_{i=1,...,n}|x_i|,x\in\mathcal{X}$. Denote $D>0$ such that $D:=\sup\{\|x-x'\|_{\infty}:x,x'\in\mathcal{X}\}$ for the diameter of $\mathcal{X}$. In this paper, we assume that $\mathcal{X}$ can be specified as $\mathcal{X}=\{x\in\mathbb{R}^n:g(x)\leq \bm{0}^{n_g}\}$ where $g:\mathbb{R}^n\rightarrow\mathbb{R}^{n_g}$ is a continuously differentiable constraint function. We have the following assumption on $g$ throughout the paper.
\begin{myassump}
    \label{assump:g_CCQ}
    Cottle Constraint Qualification (CCQ) holds at any points in $\mathcal{X}$. Namely, for any $x\in\mathcal{X}$, there is a $d\in\mathbb{R}^n$ such that  
    \begin{equation}
        \label{eq:g_CCQ}
        \nabla g(x)^\top d<\bm{0}^{n_g}
    \end{equation}
    holds.
\end{myassump}
Let $\mathscr{B}(\mathcal{X})$ be Borel $\sigma$-algebra on metric space $\mathcal{X}$. This paper uses $\mathscr{B}(\cdot)$ to denote the Borel $\sigma$-algebra on a metric space. Notice that $\left(\mathcal{X},\mathscr{B}(\mathcal{X})\right)$ is a Borel space. Let $\mu$ be a Borel probability measure on $\mathscr{B}(\mathcal{X})$. Let $M(\mathcal{X})$ be the space of Borel probability measures on metric space $\mathcal{X}$. Let $\delta$ be a random vector with support $\Delta\subseteq\mathbb{R}^s$ and $\mathbb{P}\{\cdot\}$ be the probability measurable defined on Borel $\sigma$-algebra $\mathscr{B}(\Delta)$ on $\Delta$. Let $p(\delta)$ be the probability density function associated with $\mathbb{P}\{\cdot\}$. Given a scalar function $J:\mathcal{X}\rightarrow\mathbb{R}$, and a vector-valued function $h:\mathcal{X}\times\Delta\rightarrow\mathbb{R}^m$, a chance-constrained linear program in probability measure space is formulated as:
\begin{equation} \tag{$P_{\alpha}$}
\begin{split} 
&\underset{\mu\in M(\mathcal{X})}{{\normalfont\mathsf{min}}} \,\, \int_{\mathcal{X}}J(x)\mathsf{d}\mu \\
&{\normalfont \mathsf{s.t.}}\quad  \int_{\mathcal{X}}F(x) \mathsf{d}\mu \geq 1-\alpha,
\end{split}
\end{equation} 
where $\alpha\in (0,1)$ is a given probability level and $F(x)$ is defined by
\begin{equation}
\label{eq:F_def}
    F(x):=\int_{\Delta}\mathbb{I}\{h(x,\delta)\}\mathsf{d}\mathbb{P}\{\delta\}=\int_{\Delta}\mathbb{I}\{h(x,\delta)\}p(\delta)\mathsf{d}\delta.
\end{equation}
Here, $\mathbb{I}\{y\}$ presents the indicator function written as
\begin{eqnarray*}
     \mathbb{I}\{y\}=\left\{
     \begin{array}{ll}
          1,&\text{if}\ y\leq 0,\\
          0,&\text{if}\ y>0.
          \end{array}
     \right.
\end{eqnarray*}
Note that $F(x)$ is the probability of having $h(x,\delta)\leq 0$ for given $x$. Throughout the paper, we assume the following conditions on $J(x)$ and $h(x,\delta)$.
\begin{myassump}
\label{assump:J_h}
For functions $J(x)$ and $h(x,\delta)$, the followings are supposed to be held:
\begin{itemize}
    \item[a.] $J(x)$ is continuously differentiable with respect to $x$;
    \item[b.] $h(x,\delta)$ is continuously differentiable with respect to $x$ for any $\delta\in\Delta$;
    \item[c.] For every $x\in\mathcal{X}$, $h(x,\delta)$ is continuous with respect to $\delta$;
    \item[d.] The probability density function $p(\delta)$ is continuous with respect to $\delta$;
    \item[e.] Let $\bar{h}(x,\delta):=\max_{i}h_i(x,\delta)$, $\mathsf{supp}\ p:=\mathsf{cl}\{\delta\in\Delta:p(\delta)>0\}$ ($\mathsf{cl}\{\cdot\}$ denotes the closure), and, for each $x\in\mathcal{X}$,
    \begin{equation*}
        \Delta^{\mathsf{supp}}(x):=\{\delta\in\mathsf{supp}\ p:\bar{h}(x,\delta)=0\}.
    \end{equation*}
    For each $x\in\mathcal{X}$, the following is assumed to be true:
    \begin{equation*}
        \mathbb{P}\{\Delta^{\mathsf{supp}}(x)\}=0.
    \end{equation*}
    Besides, suppose that $h(x,\delta)$ has a continuous probability density function for every $x\in\mathcal{X}$;
    \item[f.] There exists $L>0$ such that
\begin{equation*}
    \|h(x,\delta)-h(x',\delta)\|_{\infty}\leq L\|x-x'\|_{\infty},\ \forall x,x'\in\mathcal{X}\ \text{and}\ \forall\delta\in\Delta,
\end{equation*}
and
\begin{equation*}
    |J(x)-J(x')|\leq L\|x-x'\|_{\infty},\ \forall x,x'\in\mathcal{X}.
\end{equation*}
\end{itemize}
\end{myassump}
In fact, according to the content of pp. 78-79 of \cite{Kibzun}, we can obtain the continuity of $F(x)$ from Assumption \ref{assump:J_h}.

Denote the feasible region of $P_{\alpha}$ as $M_{\alpha}(\mathcal{X}):=\{\mu\in M(\mathcal{X}):\int_{\mathcal{X}}F(x) \mathsf{d}\mu \geq 1-\alpha\}$. The optimal objective function value of $P_{\alpha}$ is 
\begin{equation}
\label{eq:J_bar}
    \bar{\mathcal{J}}_{\alpha}:=\mathsf{min}\{ \int_{\mathcal{X}}J(x) \mathsf{d}\mu:\mu\in M_{\alpha}(\mathcal{X})\}.
\end{equation}
The optimal solution set of $P_{\alpha}$ is therefore written as 
\begin{equation}
\label{eq:A_alpha}
A_{\alpha}:=\{\mu\in M_{\alpha}(\mathcal{X}):\int_{\mathcal{X}}J(x) \mathsf{d}\mu = \bar{\mathcal{J}}_{\alpha}\}.    
\end{equation}
$\bar{\mu}_{\alpha}\in A_{\alpha}$ is called an optimal measure for $P_\alpha$.

\subsection{Motivation}
The motivation for addressing chance-constrained linear programs in probability measure space is from seeking an optimal stochastic policy for the optimal control problem with chance constraints, which is vital for the deployment of reliable autonomous systems by control algorithms that are robust to model misspecifications and for external disturbances \cite{Blackmore:2011,Castillo:2019,Shen_ETCI}. The optimal control problem with chance constraints aims at maximizing a reward function or minimizing a cost function with the constraints that the system state should locate in the safe area with a required probability. The deterministic policy has a fixed value in the decision domain at every time index. In contrast, the stochastic policy provides a probability measure on the decision domain at every time index. The deterministic policy can be regarded as a particular case of the stochastic policy by concentrating the probability measure on a fixed value in the decision domain. The existing techniques for addressing optimal control problems with chance constraints do not touch the essential parts of the problem and may require application-specific assumptions. For example, \cite{Castillo:2019,Hewing:2020} enforces pointwise chance constraints that ensure the independent satisfaction of each chance constraint at each time step, which leads to a more conservative solution. In general, joint chance constraints are desired, which requires all chance constraints to be satisfied jointly at all times. However, it is challenging to tackle the joint chance-constrained optimal control problem since the distribution of the state trajectory needs to be considered fully. It is possible to address the joint chance-constrained optimal control problem by using Boole's inequality \cite{Blackmore:2011,Ono:2015,Teo:2013} or performing robust optimization within the bounded model parameters obtained by specifying a confident set \cite{Lew:2022}. However, these two methods are conservative. More investigations from the viewpoint of optimization theory should be addressed to enhance new breakthroughs for optimal control with chance constraints.

Obtaining open-loop stochastic optimal policies under chance constraints can be essentially written as a chance-constrained linear program in probability measure space \cite{Thorpe:2022}. Open-loop stochastic policies mean that the stochastic policies only depend on the initial state. Unfortunately, there is still no research on solving chance-constrained linear programs in probability measure space to our knowledge. Investigating the chance-constrained linear programs in probability measure space is vital, which can give more insights into optimal control with chance constraints. 

\subsection{Related Works}

Optimization with finite chance constraints in finite-dimensional vector space is generally challenging due to the non-convexity of the feasible set and intractable reformulations \cite{Shapiro,CampiBook}. The existing research has two major streams: (1) give assumptions that the constraint functions or the distribution of random variables have some special structure, for example, linear or convex constraint functions \cite{Nemirovski}, finite sample space of random variables \cite{Luedtke:2010}, elliptically symmetric Gaussian-similar distributions \cite{Ackooij}, or (2) extract samples  \cite{Calariore:2006,Campi:2008,Luedtke:2008,Pagnoncelli:2009,Campi:2011,Teo:2017,Pena:2020,Shen:2021} or use smooth functions \cite{Geletu:2019} to approximate the chance constraints. For sample-based methods, the most famous approach in the control field is scenario approach \cite{Calariore:2006,Campi:2008,Campi:2011,Shen_ETCI,Campi_unconvex}. Scenario approach generates a deterministic optimization problem as the approximation of the original one by extracting samples from the sample space of random variables. The probability of the feasibility of the approximate solution rapidly increases to one as the sample number increases. However, the convergence of the optimality of the approximate solution is not discussed. In another sample-based method, the sample-average approach \cite{Luedtke:2008,Geletu:2019,Shen:2021,Pena:2020}, both feasibility and optimality of the approximate solution are presented. However, neither scenario approach nor sample-average approach can be directly used to solve chance-constrained linear programs in probability measure space since the deduction of the convergence of either scenario approach or sample-average approach assumes that the dimension of the decision variable must be finite. 

Optimization with chance/robust constraints in finite-dimensional vector space is also intensively studied, in which the number of chance constraints is infinite \cite{Ackooij:2019,Ackooij:2020,Chen:2021,Berthold:2022}. In \cite{Ackooij:2019}, the generalized differentiation of the probability function of infinite constraints is investigated. The optimality condition with an explicit formulation of subdifferentials is given. In \cite{Ackooij:2020}, the variational tools are applied to formulate generalized differentiation of chance/robust constraints. The method of getting the explicit outer estimations of subdifferentials from data is also established. An adaptive grid refinement algorithm is developed to solve the optimization with chance/robust constraints in \cite{Berthold:2022}. However, the above research on optimization with chance/robust constraints in finite-dimensional vector space can prove convergence only when the dimension of the decision variable is finite.

Recently, chance constraints in infinite dimensions have attracted a lot of attention. In \cite{Shaker:2018,Geletu:2020,Grandon:2022}, some essential properties, such as convexity and semi-continuity, are generalized into the chance constraints in infinite dimensions. However, the results in \cite{Shaker:2018} assume that the random variable should have a log-concave density to ensure the semicontinuity. In \cite{Grandon:2022}, the continuity of the probability function as chance constraints is proved under the assumption of continuous random distributions. The properties of chance constraints in infinite dimensions are crucial to constructing the optimality condition and implementing convergence analysis for optimization with chance constraints in infinite dimensions. In \cite{Geletu:2020}, chance-constrained optimization of elliptic partial differential equation systems is addressed by inner-outer approximation. It proves that the inner and outer approximation converges to the original problem and can provide approximate solutions with ensured convergence. However, the proof of the convergence requires the assumption that the state domain is convex. Besides, it concerns the specific problem in partial differential equation systems.

\subsection{Overview of Proposed Method and Contributions}
This paper extends the sample-based approximation method to solve chance-constrained linear programs in probability measure space. We show the relationship between chance-constrained optimization in finite-dimensional vector space and chance-constrained linear program in probability measure space. By solving a chance-constrained linear program in probability measure space, we can obtain a stochastic policy to improve the expectation of the optimal value further. We also show that the optimal objective values of the chance-constrained linear program in probability measure space and chance-constrained optimization in finite-dimensional vector space are equal if the constraints involved with random variables are required to be satisfied with probability 1. Namely, in this case, by concentrating the probability measure on an optimal solution of chance-constrained optimization in finite-dimensional vector space, we can obtain an optimal measure for the chance-constrained linear program in probability measure space. Besides, a sample approximate problem and a Gaussian mixture model approximate problem of problem $P_{\alpha}$ are proposed, by solving which the approximate solution of $P_{\alpha}$ can be obtained. The convergences of both approximate problems are investigated. Numerical examples are implemented to validate the proposed methods.

Chance-constrained linear program in probability measure space involves chance constraints in infinite dimensions. Our work differs from the \cite{Shaker:2018,Grandon:2022} in that our purpose is to provide numerical methods for solving chance-constrained linear programs in probability measure space. The properties of chance constraints in infinite dimensions are essential for convergence analysis. 

The rest of this paper is organized as follows: Section \ref{sec:main} presents two approximate problems of $P_\alpha$ and gives the main results on the convergence for each approximate problem. The proofs of the main results are presented in Section \ref{sec:proofs}. Section \ref{sec:numerical} presents the results of two numerical examples, which show the effectiveness of our proposed methods. Section \ref{sec:conclusions} concludes the whole paper.

\section{Main Results}
\label{sec:main}

This section introduces two approximate problems of $P_{\alpha}$. We also present the convergence for each approximate problem. The proofs are presented in Section \ref{sec:proofs}.

\subsection{Chance Constrained Optimization in Finite Space}
Chance-constrained optimization $Q_{\alpha}$ is an optimization problem with chance constraints in a finite-dimension vector space. The problem is written as
\begin{equation} \tag{$Q_{\alpha}$}
\begin{split}
&\underset{x\in \mathcal{X}}{{\normalfont\mathsf{min}}} \,\, J(x) \\
&{\normalfont \mathsf{s.t.}}\quad  F(x)\geq 1-\alpha,
\end{split}
\end{equation} 
where $\alpha\in (0,1)$ is a given probability level. 

Let $\mathcal{X}_{\alpha}:=\{x\in\mathcal{X}:F(x)\geq 1-\alpha\}$ be the feasible domain of $Q_{\alpha}$. Denote $\bar{J}_{\alpha}:=\mathsf{min}\{J(x):x\in\mathcal{X}_{\alpha}\}$ for the optimal objective value of $Q_{\alpha}$ and $X_{\alpha}:=\{x\in\mathcal{X}_{\alpha}:J(x)=\bar{J}_{\alpha}\}$ for the optimal solution set of $Q_{\alpha}$. We have the following assumptions over $Q_{\alpha}$ throughout the paper.
\begin{myassump}
\label{assump:x_opt}
There exists a globally optimal solution $\bar{x}$ of {\normalfont $Q_{\alpha}$} such that for any $\varepsilon>0$ there is $x\in\mathcal{X}$ such that $0<\|x-\bar{x}\|\leq\varepsilon$ and $F(x)> 1-\alpha$.
\end{myassump}

The existence of chance constraints gives rise to several difficulties. First, the structural properties of $h(x,\delta)$ might not be passed to $F(x)\geq 1-\alpha$. The feasible set $\mathcal{X}_\alpha$ can be equivalently obtained as
\begin{equation}
\mathcal{X}_\alpha=\bigcup_{\Delta_s\in\mathcal{F}}\bigcap_{\delta\in\Delta_s}\mathcal{X}_{\delta},
\end{equation}
where $\mathcal{X}_{\delta}:=\{x\in\mathcal{X}:h(x,\delta)\leq 0\}$ and $\mathcal{F}:=\{\Delta_s\in\mathscr{B}(\Delta):\mathbb{P}\mathbb\{\Delta_s\}\geq 1-\alpha\}$. Even if $h_i(x,\delta),i=1,...,m$ are all linear in $x$ for every $\delta\in\Delta$, the feasible set $\mathcal{X}_\alpha$ may not be convex due to the infinite union operations. Second, it is difficult to obtain a tractable analytical function $F(x)$ to describe the constraint or find a numerically efficient way to compute it. In most applications, $p(\delta)$ is unknown, and only samples of $\delta$ are available. We briefly review the sample-based approximation method presented in \cite{Luedtke:2008,Pagnoncelli:2009,Pena:2020}. Let $\mathcal{D}_N=\{\delta^{(1)},...,\delta^{(N)}\}$ be a set of samples randomly extracted from $\Delta$ where $N\in\mathbb{N}$. Suppose the sample extraction is independently and identically distributed. Then, $\mathcal{D}_N$ can be regarded as a random variable from the augmented sample space $\Delta^N$ with probability measure $\mathbb{P}^N\{\cdot\}$ defined on the Borel $\sigma$-algebra $\mathscr{B}(\Delta^N)$. Giving $\mathcal{D}_N$, $\epsilon\in[0,\alpha)$, and $\gamma>0$, a sample average approximate problem of $Q_{\alpha}$, defined by $\tilde{Q}_{\epsilon,\gamma}(\mathcal{D}_N)$, is written as:
\begin{equation} \tag{$\tilde{Q}_{\epsilon,\gamma}(\mathcal{D}_N)$}
\begin{split}
&\underset{x\in \mathcal{X}}{{\normalfont\mathsf{min}}} \,\, J(x) \\
&{\normalfont \mathsf{s.t.}}\quad  \frac{1}{N}\sum_{j=1}^N\mathbb{I}\{h(x,\delta^{(j)})+\gamma\}\geq 1-\epsilon. 
\end{split}
\end{equation} 
The feasible region of $\tilde{Q}_{\epsilon,\gamma}(\mathcal{D}_N)$ is defined by 
\begin{equation*}
    \tilde{\mathcal{X}}_{\epsilon,\gamma}(\mathcal{D}_N):=\{x\in\mathcal{X}:\frac{1}{N}\sum_{j=1}^N\mathbb{I}\{h(x,\delta^{(j)})+\gamma\}\geq 1-\epsilon\}.
\end{equation*}
Denote $\tilde{J}_{\epsilon,\gamma}(\mathcal{D}_N):=\mathsf{min}\{J(x):x\in\tilde{\mathcal{X}}_{\epsilon,\gamma}(\mathcal{D}_N)\}$ for the optimal objective function value of $\tilde{Q}_{\epsilon,\gamma}(\mathcal{D}_N)$ and $\tilde{X}_{\epsilon,\gamma}(\mathcal{D}_N):=\{x\in\tilde{\mathcal{X}}_{\epsilon,\gamma}(\mathcal{D}_N):J(x)=\tilde{J}_{\epsilon,\gamma}(\mathcal{D}_N)\}$ for the optimal solution set of $\tilde{Q}_{\epsilon,\gamma}(\mathcal{D}_N)$. We can regard $\tilde{J}_{\epsilon,\gamma}(\mathcal{D}_N)$ as a function $\tilde{J}_{\epsilon,\gamma}:\Delta^N\rightarrow\mathbb{R}$ for given $\epsilon$ and $\gamma$. Since $\mathcal{D}_N$ is a random variable from $\Delta^N$, $\tilde{J}_{\epsilon,\gamma}(\mathcal{D}_N)$ is consequently a random variable. The sets $\tilde{\mathcal{X}}_{\epsilon,\gamma}(\mathcal{D}_N)$ and $\tilde{X}_{\epsilon,\gamma}(\mathcal{D}_N)$ also depend on $\mathcal{D}_N$ and can be regarded as $\tilde{\mathcal{X}}_{\epsilon,\gamma}:\Delta^N\rightarrow\mathscr{B}(\mathcal{X})$ and $\tilde{X}_{\epsilon,\gamma}:\Delta^N\rightarrow\mathscr{B}(\mathcal{X})$. $\tilde{\mathcal{X}}_{\epsilon,\gamma}(\mathcal{D}_N)$ and $\tilde{X}_{\epsilon,\gamma}(\mathcal{D}_N)$ are called random sets \cite{Molchanov}. In \cite{Luedtke:2008} and \cite{Pena:2020}, the convergence analysis on $\tilde{\mathcal{X}}_{\epsilon,\gamma}(\mathcal{D}_N),\tilde{X}_{\epsilon,\gamma}(\mathcal{D}_N),\tilde{J}_{\epsilon,\gamma}(\mathcal{D}_N)$ are given. We summarize Theorem 10 of \cite{Luedtke:2008} and Theorem 3.5 of \cite{Pena:2020} as Lemma \ref{lemma:sample_average}.
\begin{mylemma}
\label{lemma:sample_average}
Suppose that Assumptions \ref{assump:J_h} and \ref{assump:x_opt} hold. Let $\epsilon\in[0,\alpha),\beta\in(0,\alpha-\epsilon)$ and $\gamma>0$. Then, 
\begin{equation*}
    {\normalfont \mathbb{P}^N}\{\tilde{\mathcal{X}}_{\epsilon,\gamma}(\mathcal{D}_N)\subseteq\mathcal{X}_{\alpha}\}\geq 1-\lceil \frac{1}{\eta}\rceil \lceil \frac{2LD}{\gamma}\rceil^{n}\exp\{-2N(\alpha-\epsilon-\beta)^2\}.
\end{equation*}
Besides, $\tilde{X}_{\epsilon,\gamma}(\mathcal{D}_N)\rightarrow X_{\alpha}$ and $\tilde{J}_{\epsilon,\gamma}(\mathcal{D}_N)\rightarrow \bar{J}_{\alpha}$ with probability 1 when $N\rightarrow\infty$, $\epsilon\rightarrow\alpha,\gamma\rightarrow 0$.
\end{mylemma}

According to Lemma \ref{lemma:sample_average}, we can obtain the solution of $Q_{\alpha}$ with probability 1 when $N\rightarrow\infty,\ \epsilon\rightarrow\alpha,\ \gamma\rightarrow 0$. A natural question arises: can we use the solution of $Q_{\alpha}$ to obtain an optimal probability measure for $P_{\alpha}$? Let $\bar{x}_{\alpha}\in X_{\alpha}$ be an optimal solution of $Q_{\alpha}$. Notice that we have $\{\bar{x}_{\alpha}\}\in\mathscr{B}(\mathcal{X})$ and thus it is possible to define a probability measure $\mu_{\bar{x}_{\alpha}}$ which satisfies that $\mu_{\bar{x}_{\alpha}}(\{\bar{x}_{\alpha}\})=\mu_{\bar{x}_{\alpha}}(\mathcal{X})=1$. Then, 
\begin{equation*}    
\int_{\mathcal{X}}J(x)\mathsf{d}\mu_{\bar{x}_{\alpha}}=\int_{\{\bar{x}_{\alpha}\}}J(x)\mathsf{d}\mu_{\bar{x}_{\alpha}}=\bar{J}_{\alpha}
\end{equation*}
and
\begin{equation*}
    \int_{\mathcal{X}}F(x) \mathsf{d}\mu_{\bar{x}_{\alpha}}=\int_{\{\bar{x}_{\alpha}\}}F(x) \mathsf{d}\mu_{\bar{x}_{\alpha}}=F(\bar{x}_{\alpha})\geq 1-\alpha.
\end{equation*}
Thus, $\mu_{\bar{x}_{\alpha}}$ is a feasible solution for $P_{\alpha}$ with objective value as $\bar{J}_{x,\alpha}$. However, $\mu_{\bar{x}_{\alpha}}$ is not sure to locate in $A_{\alpha}$. Only when $\alpha=0$, we have $\mu_{\bar{x}_{\alpha}}\in A_{\alpha}$. Notice that it is not ensured that the set $X_{\alpha}$ is a Borel measurable set. However, it is possible to find a subset $X_{\alpha}^{\text{m}}\subseteq X_{\alpha}$ that is Borel measurable. A particular example is to choose $X_{\alpha}^{\text{m}}=\{\bar{x}_{\alpha}\}$ where $\bar{x}_{\alpha}\in X_{\alpha}$ is one element in the optimal solution set. In this paper, without loss of generality, we assume that $X_{\alpha}$ is Borel measurable for all $\alpha\in[0,1]$. Besides, we also assume that $\mathcal{X}_0\neq\emptyset$. Then, $\mathcal{X}_\alpha\neq\emptyset$ holds for all $\alpha\in[0,1]$. The above content is formally summarized in Theorem \ref{theo:m1}. 
\begin{mytheo}
\label{theo:m1}
Suppose that $\mathcal{X}_{\alpha}$ is measurable for all $\alpha\in[0,1]$ and $\mathcal{X}_0\neq\emptyset$. The optimal value of problem $P_{\alpha}$ satisfies $\bar{\mathcal{J}}_{\alpha}\leq \bar{J}_{\alpha}$. Besides, if $\alpha=0$, we have
\begin{equation*}
    \bar{\mathcal{J}}_{0}=\bar{J}_{0}
\end{equation*}
and
\begin{equation}
\label{eq:mu_A_0}
    A_{0}=\{\mu\in M(\mathcal{X}):\mu(\mathcal{X})=\mu(X_0)=1\}
\end{equation}
with probability 1.
\end{mytheo}

The proof of Theorem \ref{theo:m1} is given in Section \ref{sec:proof_theo_1}. 

\begin{myrema}
Theorem \ref{theo:m1} implies that deterministic policy is optimal for robust optimal control where $\alpha=0$. 
\end{myrema}

\subsection{Sample-based Approximation} 
Let $\mathcal{X}^{\mathsf{in}}$ be the set of all interior points of $\mathcal{X}$. By using Hit-and-Run algorithm \cite{Smith} and Billiard Walk algorithm \cite{Gryazina}, uniform samples can be generated from $\mathcal{X}^{\mathsf{in}}$. For a positive integer $S\in\mathbb{N}$, let $\mathcal{C}_S:=\{x^{(1)},...,x^{(S)}\}$ be a set of uniform samples independently extracted from $\mathcal{X}^{\mathsf{in}}$. The set $\mathcal{C}_S$ is an element of the augmented space $\left(\mathcal{X}^{\mathsf{in}}\right)^S$. Since each element $x^{(i)},i=1,...,S$ in $\mathcal{C}_S$ is extracted independently, we define a $S$-fold probability $\mathbb{P}^S_{\mathsf{uni}}$ ($=\mathbb{P}_{\mathsf{uni}}\times ... \times \mathbb{P}_{\mathsf{uni}}$, $S$ times) in $\left(\mathcal{X}^{\mathsf{in}}\right)^S$. Here, $ \mathbb{P}_{\mathsf{uni}}$ is the probability measure of uniform distribution on $\mathcal{X}^{\mathsf{in}}$. 

With $\mathcal{C}_S$ and $\mathcal{D}_N$, we can obtain a sample approximate problem of $P_{\alpha}$ defined by $\tilde{P}_{\alpha}(\mathcal{C}_S,\mathcal{D}_N)$: 
\begin{equation} \tag{$\tilde{P}_{\alpha}(\mathcal{C}_S,\mathcal{D}_N)$}
\begin{split}
&\underset{\mu\in U_S}{{\normalfont\mathsf{min}}} \,\, \sum_{i=1}^{S} J(x^{(i)})\mu(i) \\
&{\normalfont \mathsf{s.t.}}\quad  \sum_{i=1}^{S}\mu(i)\frac{1}{N}\sum_{j=1}^N\mathbb{I}\{h(x^{(i)},\delta^{(j)})\}\geq 1-\alpha,
\end{split}
\end{equation}
where $U_S:=\{\mu\in\mathbb{R}^S:\sum_{i=1}^{S}\mu(i)=1,\ \mu(i)\geq 0,\ \forall i=1,...,S\}$. Define $\mathcal{F}_{\alpha}(\mathcal{C}_S,\mathcal{D}_N):=\{\mu\in U_S:\sum_{i=1}^S\mu(i)\frac{1}{N}\sum_{j=1}^N\mathbb{I}\{h(x^{(i)},\delta^{(j)})\}\geq 1-\alpha\}$ as the feasible set of $\tilde{P}_{\alpha}(\mathcal{C}_S,\mathcal{D}_N)$. Denote the optimal objective function value as 
\begin{equation*}
    \tilde{\mathcal{J}}_{\alpha}(\mathcal{C}_S,\mathcal{D}_N):=\mathsf{min}\{ \sum_{i=1}^{S} J(x^{(i)})\mu(i):\mu\in\mathcal{F}_{\alpha}(\mathcal{C}_S,\mathcal{D}_N)\}.
\end{equation*}
Denote the optimal solution set for $\tilde{P}_{\alpha}(\mathcal{C}_S,\mathcal{D}_N)$ as 
\begin{equation*}
\tilde{A}_{\alpha}(\mathcal{C}_S,\mathcal{D}_N):=\{\mu\in\mathcal{F}_{\alpha}(\mathcal{C}_S,\mathcal{D}_N):\sum_{i=1}^{S} J(x^{(i)})\mu(i)= \tilde{\mathcal{J}}_{\alpha}(\mathcal{C}_S,\mathcal{D}_N)\}.    
\end{equation*}
Let $\tilde{\mu}_{\alpha}\in \tilde{A}_{\alpha}(\mathcal{C}_S,\mathcal{D}_N)$ be an optimal measure. The optimal value $\tilde{\mathcal{J}}_{\alpha}(\mathcal{C}_S,\mathcal{D}_N)$ depends on $\mathcal{C}_S$ and $\mathcal{D}_S$, and thus it can be regarded as a function $\tilde{\mathcal{J}}_{\alpha}:\mathcal{X}^S\times\Delta^{N}\rightarrow\mathbb{R}$. Then, $\tilde{\mathcal{J}}_{\alpha}(\mathcal{C}_S,\mathcal{D}_N)$ is a random variable. Besides, $\tilde{A}_{\alpha}(\mathcal{C}_S,\mathcal{D}_N)$ is a random set. 

The deduction of the convergences of $\tilde{\mathcal{J}}_{\alpha}(\mathcal{C}_S,\mathcal{D}_N)$ and $\tilde{A}_{\alpha}(\mathcal{C}_S,\mathcal{D}_N)$ requires another assumption on $P_\alpha$. We state the assumption after a brief introduction of weak convergence.

Define a space of continuous $\mathbb{R}$-valued functions by
\begin{equation}
\label{eq:C_X_R}
    \mathscr{C}(\mathcal{X},\mathbb{R}):=\{f:\mathcal{X}\rightarrow\mathbb{R} |f\ \text{is continuous}\}.
\end{equation}
It is able to define a metric on $\mathscr{C}(\mathcal{X},\mathbb{R})$ by
\begin{equation}
\label{eq:tau}
    \tau(f,f'):=\|f-f'\|_\infty,
\end{equation}
where $\|f\|_\infty$ is defined as
\begin{equation*}
    \|f\|_\infty:=\sup_{x\in\mathcal{X}}|f(x)|.
\end{equation*}
The metric $\tau(\cdot,\cdot)$ turns $\mathscr{C}(\mathcal{X},\mathbb{R})$ into a complete metric space.

The weak convergence of probability measures is defined as follows \cite{Billingsley}.
\begin{mydef}
    \label{def:weak_convergence}
    Let $\{\mu_k\}_{k=0}^{\infty}$ be a sequence in $M(\mathcal{X})$. We say that $\{\mu_k\}_{k=0}^{\infty}$ converges weakly to $\mu$ if
    \begin{equation}
        \label{eq:weak_converge_def}
        \lim_{k\rightarrow\infty}\left|\int_{\mathcal{X}}f(x)\mathsf{d}\mu_k-\int_{\mathcal{X}}f(x)\mathsf{d}\mu\right|=0,\ \text{for all}\ f\in\mathscr{C}(\mathcal{X},\mathbb{R}).
    \end{equation}
\end{mydef}
Since $\mathcal{X}$ is compact, $M(\mathcal{X})$ can be proved to be weakly compact by Riesz representation theorem \cite{Billingsley}. Therefore, giving any sequence of $\{\mu_k\}_{k=0}^{\infty}\subset M(\mathcal{X})$, there is a subsequence which converges weakly to some $\mu\in M(\mathcal{X})$ in the sense of Definition \ref{def:weak_convergence}. By Assumption \ref{assump:J_h}, we have that $J(x)$ and $F(x)$ are continuous with respect to $x$. Therefore, if $\{\mu_k\}_{k=0}^{\infty}$ converges weakly to $\mu$, \eqref{eq:weak_converge_def} also holds for $J(x)$ or $F(x)$. We give the following assumption on Problem $P_\alpha$.
\begin{myassump}
\label{assump:P_alpha_opt_solution}
There exists a globally optimal solution $\mu^*\in A_\alpha$ of Problem $P_\alpha$ such that for any $\delta>0$ there is $\mu\in M(\mathcal{X})$ such that $\int_{\mathcal{X}}F(x)\mathsf{d}\mu>1-\alpha$ and $\mathcal{W}(\mu,\mu^*)\leq \delta$, where $\mathcal{W}(\mu,\mu^*)$ is defined by
\begin{equation}
    \label{eq:distance_mu}
    \mathcal{W}(\mu,\mu^*)=\left|\int_{\mathcal{X}}J(x)\mathsf{d}\mu-\int_{\mathcal{X}}J(x)\mathsf{d}\mu^*\right|.
\end{equation}
\end{myassump}

As $S,N\rightarrow\infty$, the convergence analysis on $\tilde{\mathcal{J}}_{\alpha}(\mathcal{C}_S,\mathcal{D}_N)$ and $\tilde{A}_{\alpha}(\mathcal{C}_S,\mathcal{D}_N)$ is summarized in Theorem \ref{theo:m2}. 
\begin{mytheo}
\label{theo:m2}
Consider Problem $P_\alpha$ with $\alpha>0$. Suppose Assumptions \ref{assump:g_CCQ}, \ref{assump:J_h}, \ref{assump:x_opt}, and \ref{assump:P_alpha_opt_solution} hold. As $S,N\rightarrow\infty$, we have
\begin{equation*}
    \liminf_{S,N\rightarrow\infty} \tilde{\mathcal{J}}_{\alpha}(\mathcal{C}_S,\mathcal{D}_N)=\bar{\mathcal{J}}_{\alpha},
\end{equation*}
with probability 1. Besides, As $S,N\rightarrow\infty$, we have $\tilde{A}_{\alpha}(\mathcal{C}_S,\mathcal{D}_N)\subset M_{\alpha}(\mathcal{X}):=\{\mu\in M(\mathcal{X}):\int_{\mathcal{X}}F(x) \mathsf{d}\mu \geq 1-\alpha\}$ with probability 1.
\end{mytheo}
The proof of Theorem \ref{theo:m2} is given in Section \ref{sec:proof_theo_2}.

\subsection{Gaussian Mixture Model-based Approximation}
Another option of approximation is to constrain the choice of $\mu$ in $M_{\theta}(\mathcal{X})\subseteq M(\mathcal{X})$. Here, $M_\theta(\mathcal{X})$ is defined as
\begin{equation*}
    M_{\theta}(\mathcal{X}):=\{ \mu\in M(\mathcal{X}): \mu(X)=\int_X p_{\theta}(x)\mathsf{d}x,\ \forall X\subseteq\mathcal{X}\},
\end{equation*}
where the probability density function $p_{\theta}(x)$ is written as
\begin{equation}
\label{eq:p_theta}
    p_{\theta}(x)=\sum_{i=1}^{L}\omega_i\phi(x,m_i,\Sigma_i).
\end{equation}
Here, $\omega_i\in[0,1],\forall i=1,..,L$, $\sum_{i=1}^{L}\omega_i=1$, and $\phi(x,m_i,\Sigma_i)$ is multivariate Gaussian distribution written by
\begin{equation*}
    \phi(x,m_i,\Sigma_i)=\frac{1}{(2\pi)^{n/2}|\Sigma_i|^{1/2}}\exp(-\frac{1}{2}(x-m_i)^\top\Sigma_i^{-1}(x-m_i)).
\end{equation*}
The notation $\theta$ denotes the parameter vector, including all the unknown parameters in $\omega_i,m_i,\Sigma_i,\forall i=1,...,L$. Denote the dimension of $\theta$ as $n_{\theta}$. The feasible domain of $\theta$ is denoted by
\begin{equation*}
    \Theta:=\{\theta\in\mathbb{R}^{n_\theta}:\sum_{i=1}^{L}\omega_i=1,\ \omega_i\geq 0\}.
\end{equation*}

Then, given a data set $\mathcal{D}_N$ and the number of Gaussian distributions $L$, we can obtain a Gaussian mixture model-based approximate problem of $P_{\alpha}$ defined by $\hat{P}_{\alpha}(L,\mathcal{D}_N)$: 
\begin{equation} \tag{$\hat{P}_{\alpha}(L,\mathcal{D}_N)$}
\begin{split}
&\underset{\theta\in\Theta}{{\normalfont\mathsf{min}}} \,\, \int_{\mathcal{X}}J(x) p_{\theta}(x)\mathsf{d}x \\
&{\normalfont \mathsf{s.t.}}\quad  \int_{\mathcal{X}}\sum_{j=1}^N\frac{1}{N}\mathbb{I}\{h(x,\delta^{(j)})\}p_{\theta}(x)\mathsf{d}x\geq 1-\alpha.
\end{split}
\end{equation} 
Denote the feasible set of $\hat{P}_{\alpha}(L,\mathcal{D}_N)$ as  
\begin{equation*}
    \Theta_{\alpha}(L,\mathcal{D}_N):=\{\theta\in\Theta:\int_{\mathcal{X}}\sum_{j=1}^N\frac{1}{N}\mathbb{I}\{h(x,\delta^{(j)})\}p_{\theta}(x)\mathsf{d}x\geq 1-\alpha\},
\end{equation*}
and the optimal objective value as
\begin{equation*}
  \hat{\mathcal{J}}_{\alpha}(L,\mathcal{D}_N):=\mathsf{min}\{\int_{\mathcal{X}}J(x) p_{\theta}(x)\mathsf{d}x:\theta\in \Theta_{\alpha}(L,\mathcal{D}_N)\}. 
\end{equation*}
Besides, the optimal solution set is
\begin{equation*}
    \hat{\Theta}_{\alpha}(L,\mathcal{D}_N):=\{\theta\in\Theta_{\alpha}(L,\mathcal{D}_N):\int_{\mathcal{X}}J(x) p_{\theta}(x)\mathsf{d}x=\hat{\mathcal{J}}_{\alpha}(L,\mathcal{D}_N)\}.
\end{equation*}
The optimal objective value $\hat{\mathcal{J}}_{\alpha}(L,\mathcal{D}_N)$ depends on the number of used Gaussian models and the data set $\mathcal{D}_N$. Since data set $\mathcal{D}_N$ is essentially random variable with support $\Delta^N$, $\hat{\mathcal{J}}_{\alpha}(L,\mathcal{D}_N)$ is also a random variable. The set $\hat{\Theta}_{\alpha}(L,\mathcal{D}_N)$ is a random set.

As $L,N\rightarrow\infty$, optimality and feasibility of using the optimal solution of $\hat{P}_{\alpha}(L,\mathcal{D}_N)$ are summarized in Theorem \ref{theo:m3}. 
\begin{mytheo}
\label{theo:m3}
Consider Problem $P_\alpha$ with $\alpha>0$. Suppose Assumptions \ref{assump:g_CCQ}, \ref{assump:J_h}, \ref{assump:x_opt}, and \ref{assump:P_alpha_opt_solution} hold. As $L,N\rightarrow\infty$, we have
    \begin{equation*}
      \liminf_{L,N\rightarrow\infty} \hat{\mathcal{J}}_{\alpha}(L,\mathcal{D}_N)=\bar{\mathcal{J}}_{\alpha},
    \end{equation*}
with probability 1. Besides, let $\hat{\theta}\in \hat{\Theta}_{\alpha}(L,\mathcal{D}_N)$ be an optimal solution of $\hat{P}_{\alpha}(L,\mathcal{D}_N)$. The corresponding probability density function is $p_{\hat{\theta}}(x)$ and the obtained probability measure is
\begin{equation*}
   \mu_{\hat{\theta}}(X):=\int_{X}p_{\hat{\theta}}(x){\normalfont\mathsf{d}}x,\ \forall X\subseteq\mathcal{X}.
\end{equation*}
We have $\mu_{\hat{\theta}}\in M_{\alpha}(\mathcal{X}):=\{\mu\in M(\mathcal{X}):\int_{\mathcal{X}}F(x) {\normalfont\mathsf{d}}\mu \geq 1-\alpha\}$ with probability 1 as $L,N\rightarrow\infty$.
\end{mytheo}
The proof of Theorem \ref{theo:m3} is given in Section \ref{sec:proof_theo_3}.

\section{Proofs of Main Results}
\label{sec:proofs}
\subsection{Proof of Theorem \ref{theo:m1}}
\label{sec:proof_theo_1}
\begin{proof}{(Theorem \ref{theo:m1}).}
Define a measure by $\bar{\mu}_{\alpha}(\cdot)$, which satisfies that $\bar{\mu}_{\alpha}(X_{\alpha})=1$. Then, we have 
\begin{align*}
    \int_{\mathcal{X}_{\alpha}} J(x)\mathsf{d}\bar{\mu}_{\alpha}=\int_{X_{\alpha}} J(x)\mathsf{d}\bar{\mu}_{\alpha}=\bar{J}_{\alpha}.
\end{align*}
Besides, for the constraint, we have
\begin{align*}
    \int_{\mathcal{X}}F(x) \mathsf{d}\bar{\mu}_{\alpha}(x)=\int_{X_{\alpha}}F(x)\mathsf{d}\bar{\mu}_{\alpha}(x)\geq 1-\alpha.
\end{align*}
Then, $\bar{\mu}_{\alpha}(\cdot)\in M_{\alpha}(\mathcal{X})$ holds. Thus, we have $\bar{\mathcal{J}}_{\alpha}\leq \int_{\mathcal{X}} J(x)\mathsf{d}\bar{\mu}_{\alpha}=\bar{J}_{\alpha}$. 

When $\alpha=0$, let $\mathcal{X}^c_{0}=\{x\in\mathcal{X}:F(x)<1\}$ be the complement set of $\mathcal{X}_0$, namely, $\mathcal{X}^c_{0}\bigcup\mathcal{X}_0=\mathcal{X}$ and $\mathcal{X}^c_{0}\bigcap\mathcal{X}_0=\emptyset$. Notice that $\mathcal{X}^c_{0}$ is Borel measurable since $\mathcal{X}_0$ is Borel measurable. Suppose that there is $\tilde{\mu}(\cdot)\in M_{0}(\mathcal{X})$ such that $\tilde{\mu}(\mathcal{X}^c_{0})>0$. Then, 
\begin{equation}
    \int_{\mathcal{X}}F(x) \mathsf{d}\tilde{\mu}(x)=\int_{\mathcal{X}_0}F(x)\mathsf{d}\tilde{\mu}(x)+\int_{\mathcal{X}_0^c}F(x)\mathsf{d}\tilde{\mu}(x)<\tilde{\mu}(\mathcal{X}_0)+\tilde{\mu}(\mathcal{X}_0^c)= 1,
\end{equation}
which conflicts with that $\tilde{\mu}\in M_0(\mathcal{X})$. Therefore, we have $\mu(\mathcal{X}^c_{0})=0$ for all $\mu\in M_0(\mathcal{X})$, which implies that $\int_{\mathcal{X}}J(x)\mathsf{d}\mu=\int_{\mathcal{X}_0}J(x)\mathsf{d}\mu$ for all $\mu\in M_0(\mathcal{X})$. 

Notice that $X_0$ is a Borel measurable set. Let $\mu^*_0(\cdot)\in A_0$ be an optimal probability measure for $P_0$ and suppose $\mu^*_{0}(X_{0})<1$ for deriving the contradiction. Thus, $\mu^*_{0}(\mathcal{X}\setminus X_{0})>0$.
The corresponding objective function is
\begin{align}
\label{eq:cor_obj_fun}
\int_{\mathcal{X}} J(x)\mathsf{d}\mu^*_{0} & = \int_{\mathcal{X}_{0}} J(x)\mathsf{d}\mu^*_{0} \\
 & = \int_{X_{0}}J(x)\mathsf{d}\mu^*_{0} +  \int_{\mathcal{X}_{0}\setminus X_{0}}J(x)\mathsf{d}\mu^*_{0}\nonumber \\
 & = \int_{X_{0}}\bar{J}_0\mathsf{d}\mu^*_{0} +  \int_{\mathcal{X}_{0}\setminus X_{0}}J(x)\mathsf{d}\mu^*_{0}\ \ \left(\because\right) \left(J(x)=\bar{J}_0,\forall x\in X_{0}\right)\nonumber \\
 & = \bar{J}_0\int_{X_{0}}\mathsf{d}\mu^*_{0} +  \int_{\mathcal{X}_{0}\setminus X_{0}}J(x)\mathsf{d}\mu^*_{0} \nonumber \\
 &=\mu^*_{0}(X_{0})\cdot \bar{J}_{0} + \int_{\mathcal{X}_{0}\setminus X_{0}}J(x)\mathsf{d}\mu^*_{0}. \nonumber
\end{align}

Denote a measure by $\bar{\mu}_{0}(\cdot)$, which satisfies that $\bar{\mu}_{0}(X_{0})=1$. Then, we have 
\begin{align}
\int_{\mathcal{X}} J(x)\mathsf{d}\bar{\mu}_{0}-\int_{\mathcal{X}} J(x)d\mu^*_{0}&=\int_{\mathcal{X}_{0}} J(x)\mathsf{d}\bar{\mu}_{0}-\int_{\mathcal{X}_{0}} J(x)d\mu^*_{0} \\
 & = \bar{J}_{0} -  \mu^*_{0}(X_{0}) \bar{J}_{0} - \int_{\mathcal{X}_{0}\setminus X_{0}}J(x)\mathsf{d}\mu^*_{0}\nonumber \\ 
 &=(1-\mu^*_{0}(X_{0}))\cdot \bar{J}_{0}-\int_{\mathcal{X}_{0}\setminus X_{0}}J(x)\mathsf{d}\mu^*_{0}\nonumber \\ 
 &=\int_{\mathcal{X}_{0}\setminus X_{0}}(\bar{J}_{0}-J(x))\mathsf{d}\mu^*_{0}\nonumber \\
 &<\int_{\mathcal{X}_{0}\setminus X_{0}}(J(x)-J(x))\mathsf{d}\mu^*_{0}=0. \nonumber
\end{align}
Thus, $\mu^*_{0}(\cdot)$ is not the optimal measure. Therefore, \eqref{eq:mu_A_0} holds, which leads to $\bar{\mathcal{J}}_0=\bar{J}_0$.\ \ \ \ \ \ \ \ \ \ \ \ \ \ \ \ \ \ \ \ \ \ \ \ \ \ \ \ \ \ \ \ \ \ \ \ \ \ \ \ \ \ \ \ \ \ \ \ \ \ \ \ \ \ \ \ \ \ \ \ \ \ \ \ \ \ \ \ \ \ \ \ \ \ \ \ \ \ \ \ \ \ \ \ \ \ \ \ \ $\blacksquare$
\end{proof}

\subsection{Proof of Theorem \ref{theo:m2}}
\label{sec:proof_theo_2}
\begin{mylemma}
\label{lemma:uniform_random_search}
Suppose that Assumption \ref{assump:g_CCQ} holds. For any $x\in\mathcal{X}$, denote a set as
\begin{equation*}
  \mathcal{B}_{\varepsilon}(x):=\{y\in\mathcal{X}:\|x-y\|\leq\varepsilon\}
\end{equation*}
where $\varepsilon>0$ is radius. For any $\varepsilon>0$, we have
\begin{equation}
\label{eq:lim_P_S}
  \lim_{S\rightarrow\infty} \mathbb{P}^S_{\mathsf{uni}}\{\mathcal{C}_S\bigcap\mathcal{B}_\varepsilon(x)\neq\emptyset\}=1.
\end{equation}
\end{mylemma}
\begin{proof}{(Lemma \ref{lemma:uniform_random_search}).}
First, we show that the interior point set $\mathcal{X}^{\mathsf{in}}$ is not empty when Assumption \ref{assump:g_CCQ} holds. Let $\bar{x}\in\mathcal{X}$ and thus we have
\begin{equation}
    \label{eq:bar_x_ine_g}
    g(\bar{x})\leq\bm{0}^{n_g}.
\end{equation}
By Assumption \ref{assump:g_CCQ}, CCQ holds at $\bar{x}$. Thus, there exists $d\in\mathbb{R}^n$ such that
\begin{equation}
    \label{eq:bar_x_ine_nable}
    \nabla g(\bar{x})^\top d<\bm{0}^{n_g}.
\end{equation}
Notice that \eqref{eq:bar_x_ine_g} and \eqref{eq:bar_x_ine_nable} directly give 
\begin{equation}
    \label{eq:bar_x_ine_g_nable}
    g(\bar{x}) + \nabla g(\bar{x})^\top d<\bm{0}^{n_g}.
\end{equation}
Since $g(\cdot)$ is continuously differentiable, there exists a small enough $\bar{\xi}>0$ such that $g(\bar{x}+\xi d)<0$ holds for any $\xi\in (0,\bar{\xi})$ and thus $\bar{x}+\xi d\in\mathcal{X}^{\mathsf{in}}$. It implies that $\mathcal{X}^{\mathsf{in}}$ is not empty.

We start from discussing $\mathbb{P}^S_{\mathsf{uni}}\{\mathcal{C}_S\bigcap\mathcal{B}_\varepsilon(x)\neq\emptyset\}$ for $x\in\mathcal{X}^{\mathsf{in}}$. Notice that $\mathcal{X}$ is compact and $\mathcal{C}_S$ is a set of uniform samples extracted from $\mathcal{X}^{\mathsf{in}}$. Thus, for any $x\in\mathcal{X}^{\mathsf{in}}$, the probability that a sample $x^{(i)}\in\mathcal{C}_S,\ i=1,..,S$ locates in $\mathcal{B}_\varepsilon(x)$ is
\begin{equation*}
\mathbb{P}_{\mathsf{uni}}\{x^{(i)}\in\mathcal{B}_\varepsilon(x)\}>0.
\end{equation*}
Then, 
\begin{align}
\label{eq:x_in_inequi_prob}
\mathbb{P}^S_{\mathsf{uni}}\{\mathcal{C}_S\bigcap\mathcal{B}_\varepsilon(x)\neq\emptyset\}&=1-\mathbb{P}^S_{\mathsf{uni}}\{\mathcal{C}_S\bigcap\mathcal{B}_\varepsilon(x)=\emptyset\} \\
&\geq 1-\left(1-\mathbb{P}_{\mathsf{uni}}\{x^{(i)}\in\mathcal{B}_\varepsilon(x)\}\right)^S. \nonumber
\end{align}
If $S\rightarrow\infty$, we have $\mathbb{P}^S_{\mathsf{uni}}\{\mathcal{C}_S\bigcap\mathcal{B}_\varepsilon(x)\neq\emptyset\}\geq 1$, which implies \eqref{eq:lim_P_S}. 

Then, we discuss $\mathbb{P}^S_{\mathsf{uni}}\{\mathcal{C}_S\bigcap\mathcal{B}_\varepsilon(x)\neq\emptyset\}$ for $x\in \partial\mathcal{X}$, where $\partial\mathcal{X}$ defines the boundary of $\mathcal{X}$. Let $x\in\partial\mathcal{X}$ be a boundary point. Again, by Assumption \ref{assump:g_CCQ}, $x$ satisfies the CCQ. By replacing $\bar{x}$ in \eqref{eq:bar_x_ine_g} and \eqref{eq:bar_x_ine_g_nable} by $x$, we have that there exists a small enough $\bar{\xi}>0$ such that $g(x+\xi d)<0$ holds for any $\xi\in (0,\bar{\xi})$ and thus $x+\xi d\in\mathcal{X}^{\mathsf{in}}$. Let $\varepsilon_1\in(0,\bar{\xi})$ and we can find $x':=x+\xi d\in\mathcal{B}_{\varepsilon_1}(x)\bigcap\mathcal{X}^{\mathsf{in}}$ with a small enough $\xi$. Besides, the probability that a sample $x^{(i)}\in\mathcal{C}_S,\ i=1,...,S$ locates in $\mathcal{B}_{\varepsilon_1}(x')$ satisfies that $\mathbb{P}_{\mathsf{uni}}\{x^{(i)}\in\mathcal{B}_{\varepsilon_1}(x')\}>0$. Thus, we have $\mathbb{P}_{\mathsf{uni}}\{x^{(i)}\in\mathcal{B}_{2\varepsilon_1}(x)\}>0$. Let $\varepsilon_1=\varepsilon/2$, and we can obtain \eqref{eq:x_in_inequi_prob} for a boundary point of $\mathcal{X}$, which completes the proof.
    \ \ \ \ \ \ \ \ \ \ \ \ \ \ \ \ \ \ \ \ \ \ \ \ \ \ \ \ \ \ \ \ \ \ \ \ \ \ \ \ \ \ \ \ \ \ \ \ \ \ \ \ \ \ \ \ \ \ \ \ \ \ \ $\blacksquare$
\end{proof}

With sample set $\mathcal{C}_{S}=\{x^{(1)},...,x^{(S)}\}$, a sample average approximate problem of $P_{\alpha}$, defined by $\breve{P}_{\alpha}(\mathcal{C}_{S})$, is written as:
\begin{equation} \tag{$\breve{P}_{\alpha}(\mathcal{C}_{S})$}
\begin{split}
&\underset{\mu\in U_S}{{\normalfont\mathsf{min}}} \,\, \sum_{i=1}^SJ(x^{(i)})\mu(i) \\
&{\normalfont \mathsf{s.t.}}\quad \sum_{i=1}^{S}\mu(i)F(x^{(i)})\geq 1-\alpha,
\end{split}
\end{equation} 
where $U_S:=\{\mu\in\mathbb{R}^S:\sum_{i=1}^{S}\mu(i)=1,\ \mu(i)\geq 0,\ \forall i=1,...,S\}$. Denote the feasible region of problem $\breve{P}_{\alpha}(\mathcal{C}_{S})$ as 
\begin{equation*}
    \breve{\mathcal{F}}_{\alpha}(\mathcal{C}_S):=\{\mu\in U_S:\sum_{i=1}^{S}\mu(i)F(x^{(i)})\geq 1-\alpha.\}.
\end{equation*}
Then, the optimal objective function value of $\breve{P}_{\alpha}(\mathcal{C}_{S})$ is defined by
\begin{equation*}
    \breve{\mathcal{J}}_{\alpha}(\mathcal{C}_S):=\mathsf{min}\{ \sum_{i=1}^SJ(x^{(i)})\mu(i):\mu\in \breve{\mathcal{F}}_{\alpha}(\mathcal{C}_S)\}.
\end{equation*}
The optimal solution set for $\breve{P}_{\alpha}(\mathcal{C}_{S})$ is therefore defined by
\begin{equation*}
    \breve{A}_{\alpha}(\mathcal{C}_S):=\{\mu\in\breve{\mathcal{F}}_{\alpha}(\mathcal{C}_S):\sum_{i=1}^SJ(x^{(i)})\mu(i)=\breve{\mathcal{J}}_{\alpha}(\mathcal{C}_S)\}.
\end{equation*}
A measure $\breve{\mu}_{\alpha}\in \breve{A}_{\alpha}(\mathcal{C}_S)$ is called an optimal measure for $\breve{P}_{\alpha}(\mathcal{C}_{S})$. 

\begin{mytheo}
\label{theo:constraint_convergence_x_sample}
For given sample sets $\mathcal{C}_{S}$ and $\mathcal{D}_{N}$, define two functions of $\mu\in U_S$ as
\begin{equation*}
    \breve{G}_{\alpha}(\mu,\mathcal{C}_{S}):=\sum_{i=1}^{S}\mu(i)F(x^{(i)})=\sum_{i=1}^{S}\mu(i)\int_{\Delta} \mathbb{I}\{h(x^{(i)},\delta)\} p(\delta)\mathsf{d}\delta,
\end{equation*}
and
\begin{equation*}
    \tilde{G}_{\alpha}(\mu,\mathcal{C}_{S},\mathcal{D}_{N}):=\sum_{i=1}^{S}\mu(i)\frac{1}{N}\sum_{j=1}^N\mathbb{I}\{h(x^{(i)},\delta^{(j)})\}.
\end{equation*}
Then, $\tilde{G}_{\alpha}(\mu,\mathcal{C}_{S},\mathcal{D}_{N})$ uniformly converges to $\breve{G}_{\alpha}(\mu,\mathcal{C}_{S})$ on $U_S$ w.p.1, i.e.,
\begin{equation*}
    \sup_{\mu\in U}\left|\tilde{G}_{\alpha}(\mu,\mathcal{C}_{S},\mathcal{D}_{N})-\breve{G}_{\alpha}(\mu,\mathcal{C}_{S})\right|\rightarrow 0,\ \text{w.p.1}\ \text{as}\ N\rightarrow\infty.
\end{equation*}
\end{mytheo}
\begin{proof}{(Theorem \ref{theo:constraint_convergence_x_sample}).}
For any given $x^{(i)}$, $\mathbb{I}\{h(x^{(i)},\delta)\}$ is a measurable function of $\delta$. According to the strong Law of Large Numbers (LLN) \cite{Bertsekas}, we have
\begin{equation*}
    \frac{1}{N}\sum_{j=1}^N\mathbb{I}\{h(x^{(i)},\delta^{(j)})\}-\mathbb{E}\{\mathbb{I}\{h(x^{(i)},\delta)\}\}\rightarrow 0,\ \text{w.p.1}\ \text{as}\ N\rightarrow\infty,
\end{equation*}
where 
\begin{equation*}
    \mathbb{E}\{\mathbb{I}\{h(x^{(i)},\delta)\}\}=\int_{\Delta} \mathbb{I}\{h(x^{(i)},\delta)\} p(\delta)\mathsf{d}\delta.
\end{equation*}
Thus, for every $\mu\in\mathcal{U}_S$, we have
\begin{align*}
    \breve{G}_{\alpha}(\mu,\mathcal{C}_{S})-\tilde{G}_{\alpha}(\mu,\mathcal{C}_{S},\mathcal{D}_{N})&=\sum_{i=1}^S \mu(i)\left(\frac{1}{N}\sum_{j=1}^N\mathbb{I}\{h(x^{(i)},\delta^{(j)})\}-\mathbb{E}\{\mathbb{I}\{h(x^{(i)},\delta)\}\}\right)\\
    &\rightarrow \sum_{i=1}^S \mu(i)\times 0=0.\ \text{w.p.1}\ \text{as}\ N\rightarrow\infty.
\end{align*}
Uniform convergence is ensured since the set $U_S$ is compact. \ \ \ \ \ \ \ \ \ \ \ \ \ \ \ \ \ \ \ \ \ \ \ \ \ \ \ \ \ \ \ \ $\blacksquare$
\end{proof}

Nextly, we show that $\tilde{\mathcal{J}}_{\alpha}(\mathcal{C}_S,\mathcal{D}_N)$ and $\tilde{A}_{\alpha}(\mathcal{C}_S,\mathcal{D}_N)$ converge to $\breve{\mathcal{J}}_{\alpha}(\mathcal{C}_S)$ and $\breve{A}_{\alpha}(\mathcal{C}_S)$, respectively, with probability 1 as $N\rightarrow\infty$.

\begin{mytheo}
\label{theo:optimal_set_convergence_x_sample}
%Suppose Assumptions \ref{assump:g_CCQ}, \ref{assump:J_h}, and \ref{assump:x_opt} hold. 
Consider Problem $P_\alpha$ with $\alpha>0$. Assume that there exists a $x^{(i)}\in\mathcal{C}_S$ that satisfies $F(x^{(i)})>1-\alpha$. As $N\rightarrow\infty$, $\tilde{\mathcal{J}}_{\alpha}(\mathcal{C}_S,\mathcal{D}_N)\rightarrow\breve{\mathcal{J}}_{\alpha}(\mathcal{C}_S)$ and $\tilde{A}_{\alpha}(\mathcal{C}_S,\mathcal{D}_N)\rightarrow \breve{A}_{\alpha}(\mathcal{C}_S)$ w.p.1. 
\end{mytheo}
\begin{proof}{(Theorem \ref{theo:optimal_set_convergence_x_sample}).}
The set $U_S$ is a compact set. The objective function $\sum_{i=1}^SJ(x^{(i)})\mu(i)$ is a linear function of $\mu\in U_S$. Besides, $F(x^{(i)})$ is a constant value within $[0,1]$ for a fixed $x^{(i)}$, which makes the constraint function $\breve{G}_{\alpha}(\mu,\mathcal{C}_{S})$ a linear function of $\mu\in U_S$. Therefore, $\breve{P}_{\alpha}(\mathcal{C}_S)$ is a linear program. Due to the assumption that there exists $x^{(i)}\in\mathcal{C}_{S}$ such that $F(x^{(i)})> 1-\alpha$, there is $\mu\in U_S$ such that $\breve{G}_{\alpha}(\mu,\mathcal{C}_{S})>1-\alpha$ and thus $\breve{A}_{\alpha}(\mathcal{C}_S)$ is nonempty. Since $\tilde{G}_{\alpha}(\mu,\mathcal{C}_{S},\mathcal{D}_{N})$ converges to $\breve{G}_{\alpha}(\mu,\mathcal{C}_{S})$ w.p.1 by Theorem \ref{theo:constraint_convergence_x_sample}, there exists $N_0$ large enough such that $\tilde{G}_{\alpha}(\mu,\mathcal{C}_{S},\mathcal{D}_{N_0})\geq 1-\alpha$ w.p.1. Because $\tilde{G}_{\alpha}(\mu,\mathcal{C}_{S},\mathcal{D}_{N_0})$ is a linear function of $\mu$ and $U_S$ is compact, the feasible set of $\tilde{P}_{\alpha}(\mathcal{C}_{S},\mathcal{D}_{N_0})$ is compact as well, and hence $\tilde{A}_{\alpha}(\mathcal{C}_{S},\mathcal{D}_{N_0})$ is nonempty w.p.1 for all $N\geq N_0$. 

Let $\{N_k\}^{\infty}_{k=1}$ be a sequence such that $N_k\rightarrow\infty$ and $N_k\geq N_0$ holds for every $k=1,...$. Let $\tilde{\mu}_k\in \tilde{A}_{\alpha}(\mathcal{C}_{S},\mathcal{D}_{N_0})$ such that $\tilde{G}_{\alpha}(\tilde{\mu}_k,\mathcal{C}_S,\mathcal{D}_{N_k})\geq 1-\alpha$, and $\sum_{i=1}^SJ(x^{(i)})\tilde{\mu}_k(i)=\tilde{\mathcal{J}}_{\alpha}(\mathcal{C}_S,\mathcal{D}_{N_k})$. Let $\tilde{\mu}$ be any cluster point of $\{\tilde{\mu}_k\}_{k=1}^{\infty}$. Let $\{\tilde{\mu}_{t}\}_{t=1}^{\infty}$ be a subsequence converging to $\tilde{\mu}$. By Theorem \ref{theo:constraint_convergence_x_sample}, we have \begin{equation*}
    \breve{G}_{\alpha}(\tilde{\mu},\mathcal{C}_S)=\lim_{t\rightarrow\infty}\tilde{G}_{\alpha}(\tilde{\mu}_t,\mathcal{C}_S,\mathcal{D}_{N_t}),\ \text{w.p.1}.
\end{equation*}
Therefore, $\breve{G}_{\alpha}(\tilde{\mu},\mathcal{C}_S)\geq 1-\alpha$ and $\tilde{\mu}$ is feasible for problem $\breve{P}_{\alpha}(\mathcal{C}_S)$ which implies $\sum_{i=1}^SJ(x^{(i)})\tilde{\mu}(i)\geq \breve{\mathcal{J}}_{\alpha}(\mathcal{C}_S)$. Note that $\tilde{\mu}_t\rightarrow \tilde{\mu}$ w.p.1, which implies that
\begin{equation*}
    \lim_{t\rightarrow\infty} \tilde{\mathcal{J}}_\alpha(\mathcal{C}_S,\mathcal{D}_{N_t})=\lim_{t\rightarrow\infty} \sum_{i=1}^SJ(x^{(i)})\tilde{\mu}_t(i)=\sum_{i=1}^SJ(x^{(i)})\tilde{\mu}(i)\geq \breve{\mathcal{J}}_{\alpha}(\mathcal{C}_S),\ \text{w.p.1}.
\end{equation*}
Since this is true for an arbitrary point of $\{\tilde{\mu}_k\}_{k=1}^{\infty}$ in the compact set $\mathcal{U}_S$, we have
\begin{equation}
\label{eq:conver_opt_obj_1}
    \lim_{k\rightarrow\infty} \tilde{\mathcal{J}}_\alpha(\mathcal{C}_S,\mathcal{D}_{N_k})=\lim_{k\rightarrow\infty} \sum_{i=1}^SJ(x^{(i)})\tilde{\mu}_k(i)\geq \breve{\mathcal{J}}_{\alpha}(\mathcal{C}_S),\ \text{w.p.1}.
\end{equation}

Besides, we know that there exists a globally optimal solution of $\breve{P}_{\alpha}(\mathcal{C}_S)$, $\mu^*$, such that for any $\varepsilon>0$ there is $\mu\in U$ such that $0<\|\mu-\mu^*\|\leq\varepsilon$ and $\breve{G}_{\alpha}(\mu,\mathcal{C}_S)>1-\alpha$. Namely, there exists a sequence $\{\tilde{\mu}_t\}_{t=1}^{\infty}\subseteq U$ that converges to an optimal solution $\mu^*$ such that $\breve{G}_{\alpha}(\tilde{\mu}_t,\mathcal{C}_S)>1-\alpha$ for all $t\in\mathbb{N}$. Notice that $\tilde{G}_{\alpha}(\tilde{\mu}_t,\mathcal{C}_S,\mathcal{D}_{N_k})$ converges to $\breve{G}_{\alpha}(\tilde{\mu}_t,\mathcal{C}_S)$ w.p.1. Then, for any fixed $t$, $\exists K(t)$ such that $\tilde{G}_{\alpha}(\tilde{\mu}_t,\mathcal{C}_S,\mathcal{D}_{N_k})\geq 1-\alpha$ for every $k\geq K(t)$ w.p.1. We can assume $K(t)<K(t+1)$ for every $t$ and define the sequence $\{\tilde{\mu}_k\}_{k=K(1)}^{\infty}$ by setting $\tilde{\mu}_k=\tilde{\mu}_t$ for all $k$ and $t$ with $K(t)\leq k<K(t+1)$. Then, $\tilde{G}_{\alpha}(\hat{\mu}_k,\mathcal{C}_S,\mathcal{D}_{N_k})\geq 1-\alpha$, which implies $\tilde{\mathcal{J}}_{\alpha}(\mathcal{C}_S,\mathcal{D}_{N_k})\leq \sum_{i=1}^SJ(x^{(i)})\tilde{\mu}_k(i)$ for all $k\geq K(1)$. Thus, we have that
\begin{equation}
\label{eq:conver_opt_obj_2}
    \lim_{k\rightarrow\infty} \tilde{\mathcal{J}}_{\alpha}(\mathcal{C}_S,\mathcal{D}_{N_k})\leq \sum_{i=1}^SJ(x^{(i)})\mu^*(i)=\breve{\mathcal{J}}_{\alpha}(\mathcal{C}_S),\ \text{w.p.1}.
\end{equation}
With \eqref{eq:conver_opt_obj_1} and \eqref{eq:conver_opt_obj_2}, we conclude that $\tilde{\mathcal{J}}_{\alpha}(\mathcal{C}_S,\mathcal{D}_{N})\rightarrow\breve{\mathcal{J}}_{\alpha}(\mathcal{C}_S)$ w.p.1 as $N\rightarrow\infty$.

The proof of $\tilde{A}_{\alpha}(\mathcal{C}_S,\mathcal{D}_N)\rightarrow \breve{A}_{\alpha}(\mathcal{C}_S)$ can be referred to Theorem 5.3 of \cite{Shapiro}.\\

\ \ \ \ \ \ \ \ \ \ \ \ \ \ \ \ \ \ \ \ \ \ \ \ \ \ \ \ \ \ \ \ \ \ \ \ \ \ \ \ \ \ \ \ \ \ \ \ \ \ \ \ \ \ \ \ \ \ \ \ \ \ \ \ \ \ \ \ \ \ \ \ \ \ \ \ \ \ \ \ \ \ \ \ \ \ \ \ \ \ \ \ \ \ \ \ \ \ \ \ \ \ \ \ \ \ \ \ \ \ \ \ \ \   $\blacksquare$

\end{proof}

Nextly, we show that $\breve{\mathcal{J}}_{\alpha}(\mathcal{C}_S)$ converges to $\bar{\mathcal{J}}_{\alpha}$ with probability 1 as $S$ increases.
\begin{mytheo}
\label{theo:measure_convergence}
Suppose Assumption \ref{assump:J_h} and \ref{assump:P_alpha_opt_solution} hold. As $S\rightarrow\infty$, with probability 1, we have
\begin{equation}
\label{eq:J_converge_measure}
\liminf_{S\rightarrow\infty} \breve{\mathcal{J}}_{\alpha}(\mathcal{C}_S)=\bar{\mathcal{J}}_{\alpha}.
\end{equation}
\end{mytheo}
\begin{proof}{(Theorem \ref{theo:measure_convergence}).}
The outline of the proof of Theorem \ref{theo:measure_convergence} is summarized as follows:
\begin{itemize}
    \item[A.] Prove that the limit of lower bound of $\breve{\mathcal{J}}_{\alpha}(\mathcal{C}_S)$ is larger than $\bar{\mathcal{J}}_{\alpha}$ by \eqref{eq:J_CS_lim_J_opt_1};
    \item[B.] Prove that the limit of upper bound of $\breve{\mathcal{J}}_{\alpha}(\mathcal{C}_S)$ is smaller than $\bar{\mathcal{J}}_{\alpha}$ by \eqref{eq:J_CS_lim_J_opt_2}; 
    \begin{itemize}
        \item[B1.] Find a sequence $\{\mu_k\}_{k=1}^{\infty}$ converges weakly to an optimal solution $\mu^*$ of $P_\alpha$;
        \item[B2.] Show that $\int_{\mathcal{X}}F(x)\mathsf{d}\mu_k(x)$ and $\int_{\mathcal{X}}J(x)\mathsf{d}\mu_k(x)$ can be approximated by using discrete probability measure on $\mathcal{C}_S$, which refers to \eqref{eq:P_int_sum_approx} and \eqref{eq:J_int_sum_approx};
        \item[B3.] Show that optimal discrete probability measure on $\mathcal{C}_S$ for $\breve{P}_\alpha(\mathcal{C}_S)$ has a smaller objective value than the discrete probability measure for approximating any $\mu_k$ in B2. Then, we obtain \eqref{eq:J_CS_lim_J_opt_2}. 
    \end{itemize}
\end{itemize}

Then, we give the details of the proof. 

For any discrete probability measure $\mu^S\in\breve{\mathcal{F}}_{\alpha}(\mathcal{C}_S)$, we have
\begin{equation*}
    \int_{\mathcal{X}}F(x) \mathsf{d}\mu^S(x)=\sum_{i=1}^{S}\mu^S(\{x^{(i)}\})F(x^{(i)})\geq 1-\alpha.
\end{equation*}
Thus, $\mu^S\in M_{\alpha}(x)$. Then, it holds that
\begin{equation*}
    \sum_{i=1}^S J(x^{(i)})\mu^S(\{x^{(i)}\})=\int_{\mathcal{X}}J(x)\mathsf{d}\mu^S(x)\geq \bar{\mathcal{J}}_{\alpha},\ \forall \mu^S\in\breve{\mathcal{F}}_{\alpha}(\mathcal{C}_S).
\end{equation*}
Furthermore, with probability 1, we have
\begin{equation}
    \label{eq:J_CS_lim_J_opt_1}
    \liminf_{S\rightarrow\infty} \breve{\mathcal{J}}_\alpha(\mathcal{C}_S)\geq\bar{\mathcal{J}}_{\alpha}.
\end{equation}
Assumption \ref{assump:P_alpha_opt_solution} implies that there exists a sequence $\{\mu_k\}_{k=1}^{\infty}\subseteq M(\mathcal{X})$ that converges weakly to an optimal solution $\mu^*$ such that
\begin{equation}
\label{eq:P_mu_k}
    \int_{\mathcal{X}}F(x) \mathsf{d}\mu_k(x)> 1-\alpha
\end{equation}
for all $k\in\mathbb{N}$. Since $\{\mu_k\}_{k=1}^{\infty}$ converges weakly to $\mu^*$, we have 
\begin{equation}
    \lim_{k\rightarrow\infty}\int_{\mathcal{X}}J(x)\mathsf{d}\mu_k(x)-\int_{\mathcal{X}}J(x)\mathsf{d}\mu^*(x)=\lim_{k\rightarrow\infty}\mathcal{W}(\mu_k,\mu^*)=0.
\end{equation}
Notice that $\bar{\mathcal{J}}_\alpha=\int_{\mathcal{X}}J(x)\mathsf{d}\mu^*(x)$ by \eqref{eq:J_bar}.

For any given $\varepsilon_J>0$, $\exists K(\varepsilon_J)$, if $k\geq K(\varepsilon_J)$, 
\begin{equation*}
    \int_{\mathcal{X}}J(x)\mathsf{d}\mu_k(x)-\bar{\mathcal{J}}_{\alpha}\leq\varepsilon_J.
\end{equation*}

Let $\tilde{\mathcal{C}}_{\tilde{S}}^k:=\{\tilde{x}^{(1)}_k,...,\tilde{x}^{(\tilde{S})}_k\}$ be a sample set obtained by sampling from $\mathcal{X}$ according to probability measure $\mu_k$. By Law of Large Numbers (p. 457 of \cite{Shapiro}), for any $f\in\mathscr{C}(\mathcal{X},\mathbb{R})$, as $\tilde{S}\rightarrow\infty$, with probability 1, we have 
\begin{equation}
    \label{eq:LLN}
    \frac{1}{\tilde{S}}\sum_{i=1}^{\tilde{S}} f(\tilde{x}^{(i)}_k)\rightarrow \mathbb{E}_{x\sim\mu_k}\left\{f(x)\right\}=\int_{\mathcal{X}}f(x)\mathsf{d}\mu_k(x).
\end{equation}
Since $J(\cdot)$ and $F(\cdot)$ are also elements in $\mathscr{C}(\mathcal{X},\mathbb{R})$, \eqref{eq:LLN} also holds by replacing $f(\cdot)$ by either $J(\cdot)$ or $F(\cdot)$. Namely, for any $\tilde{\varepsilon}_1$, there exists $\tilde{S}_{\mathsf{l}}(\tilde{\varepsilon}_J)$ such that, if $\tilde{S}\geq\tilde{S}_{\mathsf{l}}(\tilde{\varepsilon}_J)$, with probability 1, the followings hold:
\begin{equation}
    \label{eq:P_int_approx}
    \left|\frac{1}{\tilde{S}}\sum_{i=1}^{\tilde{S}} F(\tilde{x}^{(i)}_k)-\int_{\mathcal{X}}F(x)\mathsf{d}\mu_k(x)\right|\leq \tilde{\varepsilon}_1,
\end{equation}
\begin{equation}
    \label{eq:J_int_approx}
    \left|\frac{1}{\tilde{S}}\sum_{i=1}^{\tilde{S}} J(\tilde{x}^{(i)}_k)-\int_{\mathcal{X}}J(x)\mathsf{d}\mu_k(x)\right|\leq \tilde{\varepsilon}_1.
\end{equation}

On the other hand, according to Lemma \ref{lemma:uniform_random_search}, as $S\rightarrow\infty$, for any $\tilde{s}\in\{1,...,\tilde{S}\}$ and $\tilde{\varepsilon}_r>0$, with probability 1, there exists a sample $x^{(i_{\tilde{s}})}\in\mathcal{C}_S:=\{x^{(1)},...,x^{(S)}\}$ such that
\begin{equation}
    x^{(i_{\tilde{s}})}\in\mathcal{B}_{\tilde{\varepsilon}_r}(\tilde{x}_{k}^{\tilde{s}}).
\end{equation}
With a little abuse of notation, let $x^{(i_{\tilde{s}})}$ be the closest sample to $\tilde{x}^{(\tilde{s})}_k$, namely, $x^{(i_{\tilde{s}})}\in\arg\min\{\|x^{(i)}-\tilde{x}^{(\tilde{s})}_k\|:x^{(i)}\in\mathcal{C}_S\}$. Define a set $I_{\tilde{S}}:=\{i_1,...,i_{\tilde{S}}\}$ as the set of index corresponding to $x^{(i_{\tilde{s}})}$. Without loss of generality, we assume that $x^{(i_{\tilde{s}})}\neq x^{(j_{\tilde{s}})}$ if $i_{\tilde{s}}\neq j_{\tilde{s}},\ i_{\tilde{s}},j_{\tilde{s}}\in I_{\tilde{S}}$. The intuitive explanation of the relationship between $\mathcal{C}_S$ and $\tilde{\mathcal{C}}^k_{\tilde{S}}$ is illustrated in Figure \ref{fig:outline_proof_theo_measure_convergence}.
\begin{figure}
\centering
\includegraphics[scale=0.375]{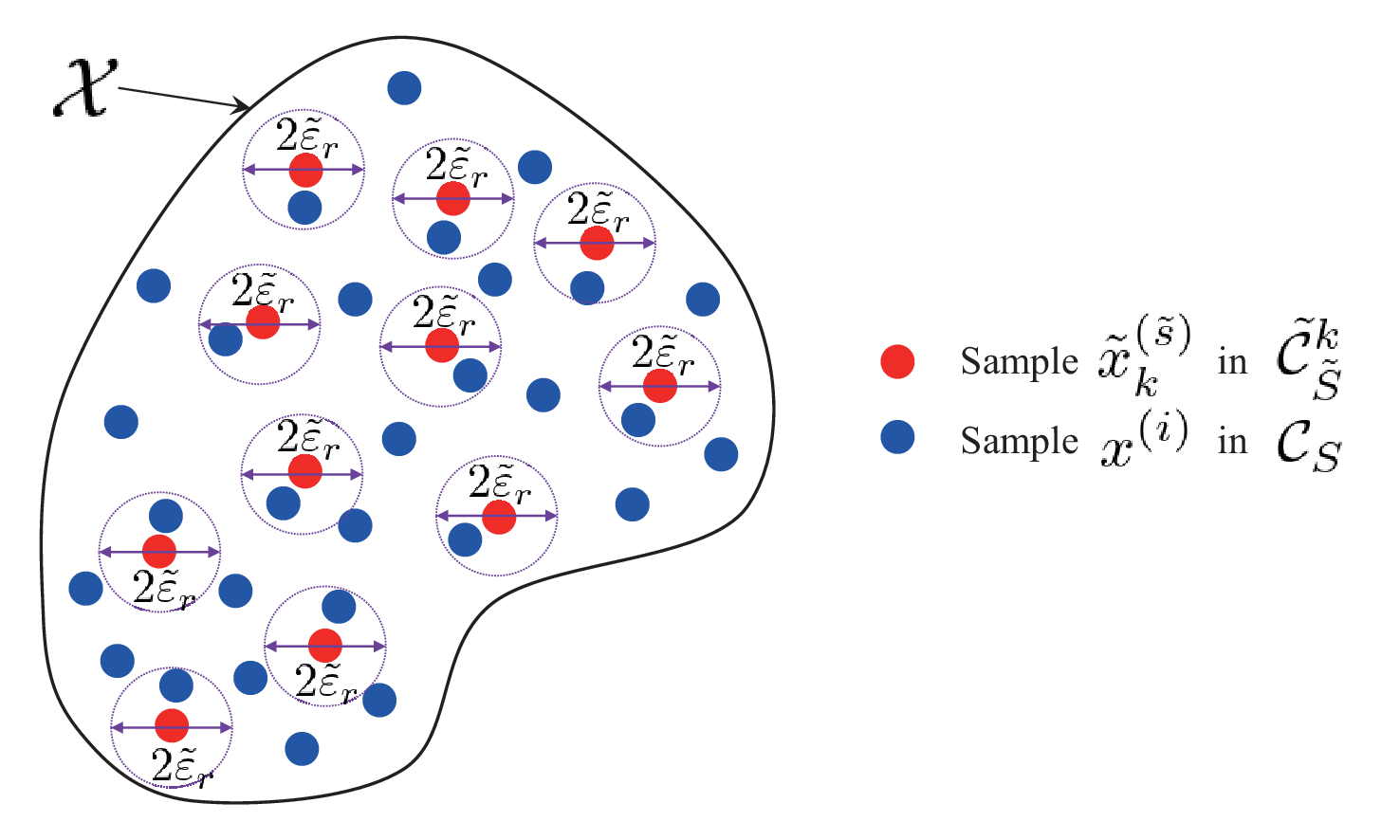}
\centering
\caption{The intuitive explanation of the relationship between $\mathcal{C}_S$ and $\tilde{\mathcal{C}}^k_{\tilde{S}}$}
\label{fig:outline_proof_theo_measure_convergence}
\end{figure}

Define a discrete probability measure $\mu^{S}_k\in \mathbb{R}^S$ such that
\begin{equation}
    \mu^S_k(i)=\frac{1}{\tilde{S}},\ \forall i\in I_{\tilde{S}},
\end{equation}
\begin{equation}
    \mu^S_k(i)=0,\ \forall i\notin I_{\tilde{S}}.
\end{equation}

For any given positive integer $\tilde{S}$ and positive number $\tilde{\varepsilon}_2$, due to the continuity of $J(\cdot)$ and $F(\cdot)$, there exists $\bar{S}_{\mathsf{l}}(\tilde{S},\tilde{\varepsilon}_2)$ such that, if $S>S_{\mathsf{l}}(\tilde{S},\tilde{\varepsilon}_2)$, with probability 1, the followings hold:
\begin{equation}
    \label{eq:P_sum_approx}
    \left|\sum_{i=1}^{S} \mu^S_k(i)F(x^{(i)})-\frac{1}{\tilde{S}}\sum_{i=1}^{\tilde{S}} F(\tilde{x}^{(i)}_k)\right|\leq \tilde{\varepsilon}_2,
\end{equation}
\begin{equation}
    \label{eq:J_sum_approx}
    \left|\sum_{i=1}^{S} \mu^S_k(i)J(x^{(i)})-\frac{1}{\tilde{S}}\sum_{i=1}^{\tilde{S}} J(\tilde{x}^{(i)}_k)\right|\leq \tilde{\varepsilon}_2.
\end{equation}
By combining \eqref{eq:P_int_approx} with \eqref{eq:P_sum_approx} and combining \eqref{eq:J_int_approx} with \eqref{eq:J_sum_approx}, then, for given $\tilde{\varepsilon}_1,\tilde{\varepsilon}_2$, there exists $\tilde{S}_{\mathsf{l}}(\tilde{\varepsilon}_1)$ and $S_{\mathsf{l}}(\tilde{S},\tilde{\varepsilon}_2)$ such that, if $\tilde{S}>\tilde{S}_{\mathsf{l}}(\tilde{\varepsilon}_1)$ and $S>S_{\mathsf{l}}(\tilde{S},\tilde{\varepsilon}_2)$, with probability 1, the following holds:
\begin{equation}
    \label{eq:P_int_sum_approx}
    \left|\sum_{i=1}^{S} \mu^S_k(i)F(x^{(i)})-\int_{\mathcal{X}}F(x)\mathsf{d}\mu_k(x)\right|\leq \tilde{\varepsilon}_1+\tilde{\varepsilon}_2,
\end{equation}
\begin{equation}
    \label{eq:J_int_sum_approx}
    \left|\sum_{i=1}^{S} \mu^S_k(i)J(x^{(i)})-\int_{\mathcal{X}}J(x)\mathsf{d}\mu_k(x)\right|\leq \tilde{\varepsilon}_1+\tilde{\varepsilon}_2.
\end{equation}
According to \eqref{eq:P_mu_k} and \eqref{eq:P_int_sum_approx}, we can find $\tilde{S}_{\mathsf{l}}(\tilde{\varepsilon}_1)$ and $S_{\mathsf{l}}(\tilde{S},\tilde{\varepsilon}_2)$ such that, if $\tilde{S}>\tilde{S}_{\mathsf{l}}(\tilde{\varepsilon}_1)$ and $S>S_{\mathsf{l}}(\tilde{S},\tilde{\varepsilon}_2)$, with probability 1, the following holds
\begin{equation}
    \label{eq:P_int_sum_approx_ineq}
    \sum_{i=1}^{S} \mu^S_k(i)F(x^{(i)})\geq 1-\alpha.
\end{equation}
Thus, $\mu^S_k$ is a feasible solution of Problem $\tilde{P}_\alpha(\mathcal{C}_S)$ and thus 
\begin{equation}
    \label{eq:mu_S_k_cost_com}
    \sum_{i=1}^{S} \mu^S_k(i)J(x^{(i)})\geq \breve{\mathcal{J}}_{\alpha}(\mathcal{C}_S).
\end{equation}
Since $\int_{\mathcal{X}}J(x)\mathsf{d}\mu_k(x)$ converges to $\bar{\mathcal{J}}_\alpha$ w.p.1 as $k\rightarrow\infty$, thus, considering \eqref{eq:J_int_sum_approx} and \eqref{eq:mu_S_k_cost_com}, we have
\begin{equation}
    \label{eq:J_CS_lim_J_opt_2}
    \limsup_{S\rightarrow\infty} \breve{\mathcal{J}}_\alpha(\mathcal{C}_S)\leq\bar{\mathcal{J}}_{\alpha}.
\end{equation}

With \eqref{eq:J_CS_lim_J_opt_1} and \eqref{eq:J_CS_lim_J_opt_2}, we have \eqref{eq:J_converge_measure}.

\ \ \ \ \ \ \ \ \ \ \ \ \ \ \ \ \ \ \ \ \ \ \ \ \ \ \ \ \ \ \ \  \ \ \ \ \ \ \ \ \ \ \ \ \ \ \ \ \ \ \ \ \ \ \ \ \ \ \ \ \ \ \ \  $\blacksquare$
\end{proof}

The proof of Theorem \ref{theo:m2} can be obtained immediately by using the results of Theorem \ref{theo:optimal_set_convergence_x_sample} and Theorem \ref{theo:measure_convergence}, which is omitted here.

\subsection{Proof of Theorem \ref{theo:m3}}
\label{sec:proof_theo_3}
Main results of \cite{Zeevi:1997} are summarized as:
\begin{mylemma}
\label{lemma:Zeevi}
Let $\mathcal{X}^+$ be a compact set. Let $p:\mathbb{R}^n\rightarrow\mathbb{R}$ be a probability density function on the domain $\mathbb{R}^n$. If there exists a positive number $\rho'>0$ such that $p\in\{p:p(x)\geq\rho',\forall x\in\mathcal{X}^+\}$, then there exists $p_{\theta}(x)$ defined by \eqref{eq:p_theta} such that
\begin{equation*}
    \lim_{L\rightarrow\infty}   \int_{\mathcal{X}^+}\left(p(x)-p_{\theta}(x)\right)^{2}\mathsf{d}x =0,
\end{equation*}
where the positive integer $L$ is the number of Gaussian kernels in \eqref{eq:p_theta}. 
\end{mylemma}

\begin{proof}{(Theorem \ref{theo:m3}).}
For given $\mathcal{C}_S,\ \mathcal{D}_{N}$ and $L$, we have problems $\tilde{P}_{\alpha}(\mathcal{C}_S,\mathcal{D}_{N})$ and $\hat{P}_{\alpha}(L,\mathcal{D}_N)$. Let $\mathcal{X}_{p,i},i=1,...,S$ be the partitions of $\mathcal{X}$, which satisfy
\begin{itemize}
    \item[(a)] $x^{(i)}\in\mathcal{X}_{p,i}$;
    \item[(b)] $\bigcup_{i=1}^{S} \mathcal{X}_{p,i}=\mathcal{X}$;
    \item[(c)] $\mathcal{X}_{p,i}\bigcap\mathcal{X}_{p,i'}=\emptyset$ with probability 1 if $i\neq i'$.
\end{itemize}
For any $\mu^{S}\in U$, we can correspondingly define a Dirac measure on $\mathcal{X}$ as
\begin{equation*}
    \mu^{S}_{\mathsf{d}}(x)=\mu^{S}(x^{(i)})\ \text{if}\ x\in\mathcal{X}_{p,i}.
\end{equation*}
Define a set of index as $I^+=\{i:\mu^{S}(x^{(i)})>0\}$. Then, we can define a compact set 
\begin{equation*}
    \mathcal{X}^+=\bigcup_{i\in I^+} \mathcal{X}_{p,i}.
\end{equation*}
According to Lemma \ref{lemma:Zeevi}, there exists a sequence $\left\{p_{\theta}(x)\right\}_L$ such that
\begin{equation*}
    \lim_{L\rightarrow\infty}   \int_{\mathcal{X}^+}\left(\mu^{S}_{\mathsf{d}}(x)-p_{\theta}(x)\right)^{2}\mathsf{d}x =0.
\end{equation*}
Thus, we have 
\begin{equation*}
    \lim_{L\rightarrow\infty}\int_{\mathcal{X}}J(x) \mathsf{d}p_{\theta}(x)=\int_{\mathcal{X}}J(x)\mathsf{d}\mu^{S}_{\mathsf{d}}(x)
\end{equation*}
and 
\begin{equation*}
    \lim_{L\rightarrow\infty}\int_{\mathcal{X}}\frac{1}{N}\sum_{j=1}^N\mathbb{I}\{h(x,\delta^{(j)})\leq 0\}p_{\theta}(x)\mathsf{d}x=\int_{\mathcal{X}}\frac{1}{N}\sum_{j=1}^N\mathbb{I}\{h(x,\delta^{(j)})\leq 0\}\mu^{S}_{\mathsf{d}}(x)\mathsf{d}x.
\end{equation*}
For any $S$ and $N$, by applying Lemma \ref{lemma:Zeevi}, we can find a sequence $\left\{p^*_{\theta}(x)\right\}_L$ such that
\begin{equation}
\label{eq:lim_L_J}
    \lim_{L\rightarrow\infty} \int_{\mathcal{X}}J(x) \mathsf{d}p_{\theta}^{*}(x)=\tilde{\mathcal{J}}_{\alpha}(\mathcal{C}_S,\mathcal{D}_N)
\end{equation}
and 
\begin{equation}
\label{eq:lim_L_F}
    \lim_{L\rightarrow\infty}\int_{\mathcal{X}}\sum_{j=1}^{N}\frac{1}{N}\mathbb{I}\{h(x,\delta^{(j)})\}p_{\theta}^{*}(x)\mathsf{d}x=\sum_{i=1}^{S}\mu^{S}(x^{(i)})\sum_{j=1}^{N}\frac{1}{N}\mathbb{I}\{h(x^{(i)},\delta^{(j)})\}\geq 1-\alpha.
\end{equation}

There exists the limit of $\hat{P}_{\alpha}(L,\mathcal{D}_N)$ that converges to $\tilde{P}_{\alpha}(\mathcal{C}_S,\mathcal{D}_N)$ as $L\rightarrow\infty$. Theorem \ref{theo:m3} can be obtained by using Theorem \ref{theo:m2}. One point should be clarified here. In Theorem \ref{theo:m2}, the convergence holds for $S\rightarrow\infty$. In Theorem \ref{theo:m3}, $L\rightarrow\infty$ is used instead since we have \eqref{eq:lim_L_J} and \eqref{eq:lim_L_F} for any $S$ increasing to infinite.\ \ \ \ \ \ \ \ \ \ \ \ \ \ \ \ \ \ \ \ \ \ \ \ \ \ \ \ \ \ \ \ \ \ \ \ \ \ \ \ \ \ \ \ \ \ \ \ \ \ \ \ \ \ \ \ \ \ \ \ \ \ \ \ \ \ \ \ \ \ \ \ \ \ \ \ \ \ \ \ \ \ \ \ \ \ \ \ \ \ \ \ \ \ \ \ \ \ \ \ \ \ \ \ \ \ \ \ \ \ \ \ \ \ \ \ \ \ \ \ \ \ \ \ $\blacksquare$
\end{proof}

\section{Numerical Examples}
\label{sec:numerical}
This section provides the results of two numerical examples to validate our proposed methods. All computations were executed on Windows 10 with 32GB RAM and an Intel(R) Core(TM) i7-1065G7 CPU running at 1.30GHz. The algorithm and all computations were implemented in Matlab R2021b. We check the performance of the following methods:
\begin{enumerate}
    \item \textbf{Dirac-Delta}: solving sample average approximate problem $\tilde{Q}_{\epsilon,\gamma}(\mathcal{D}_N)$ of $Q_\alpha$;
    \item \textbf{Sample}: solving sample-based approximate problem $\tilde{P}_\alpha(\mathcal{C}_S,\mathcal{D}_N)$ of $P_\alpha$;
    \item \textbf{GMM}: GMM-based approximate problem $\hat{P}_\alpha(L,\mathcal{D}_N)$.
\end{enumerate}
We use the terminology \textbf{Dirac-Delta} for the method of solving sample average approximate problem $\tilde{Q}_{\epsilon,\gamma}(\mathcal{D}_N)$ of $Q_\alpha$ since it equivalently gives the measure constrained to be a dirac-delta, namely, the measure is concentrated on one fixed solution.
\begin{figure}
\centering
\includegraphics[scale=0.95]{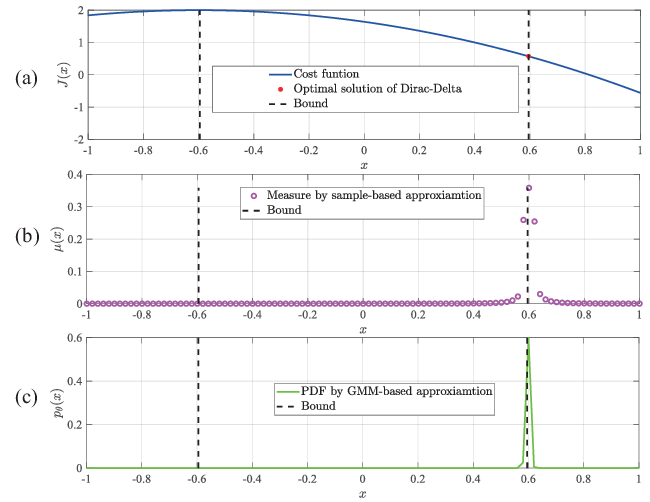}
\centering
\caption{Results of the numerical example 1: (a) profile of $J(x)$ and optimal solution obtained by \textbf{Dirac-Delta}; (b) optimal measure by \textbf{Sample}; (c) optimal probability density function obtained by \textbf{GMM}.}
\label{fig:Result_example_one}
\end{figure}
\subsection{One Dimension Example}

\begin{table}[H]
 \centering
 \label{tb:statistics} 
 \caption{Statistics of CPU time for one dimension example.}  
 \begin{tabular}{ c||c|c|c|c|c }
 \hline
 \hline
 \multicolumn{5}{c}{\textbf{Dirac-Delta}} \\
 \hline
  Size $N$ & $1000$ & $2000$ & $5000$ & $10000$\\
 \hline
 \hline
 Avg. (s)   &  $0.0105$ & $0.0136$ & $0.0154$ & $0.0169$ \\
 \hline
 Max. (s)   &  $0.0281$ & $0.0337$ & $0.0402$ & $0.0583$  \\
 \hline
 \hline
\end{tabular}
 \begin{tabular}{ c||c|c|c|c|c|c|c|c|c }
 \hline
 \hline
 \multicolumn{10}{c}{\textbf{Sample}} \\
 \hline
Size $N$ &   \multicolumn{3}{c|}{$2000$} & \multicolumn{3}{c|}{$5000$} & \multicolumn{3}{c}{$10000$}\\
 \hline
 Size $S$ & $50$ & $100$ & $200$ & $50$ & $100$ & $200$ & $50$ & $100$ & $200$  \\
 \hline
 \hline
 Avg. (s)     & $0.0113$ & $0.0123$ & $0.0140$ & $0.0119$ & $0.0145$ & $0.0185$ & $0.0143$ & $0.0162$ & $0.0229$  \\
\hline
 Max. (s)     & $0.0263$ & $0.0349$ & $0.0415$ & $0.0238$& $0.0341$ & $0.0445$ & $0.0343$ & $0.0409$ & $0.0530$ \\
 \hline
 \hline
\end{tabular}
 \begin{tabular}{ c||c|c|c|c|c|c|c|c|c }
 \hline
 \hline
 \multicolumn{10}{c}{\textbf{GMM}} \\
 \hline
Size $N$ &   \multicolumn{3}{c|}{$2000$} & \multicolumn{3}{c|}{$5000$} & \multicolumn{3}{c}{$10000$}\\
 \hline
 Size $L$ & $2$ & $4$ & $6$ & $2$ & $4$ & $6$ & $2$ & $4$ & $6$  \\
 \hline
 \hline
 Avg. (s)     & $0.0325$ & $0.0431$ & $0.0991$ & $0.0472$ & $0.0783$ & $0.1671$ & $0.0773$ & $0.1429$ & $0.2134$  \\
\hline
 Max. (s)     & $0.0721$ & $0.0784$ & $0.1562$ & $0.0793$& $0.1318$ & $0.2011$ & $0.1267$ & $0.1892$ & $0.2690$ \\
 \hline
 \hline
\end{tabular}
\end{table}
In the first numerical, we use an extremely simple example to demonstrate the concepts of Theorem \ref{theo:m1}, Theorem \ref{theo:m2}, and Theorem \ref{theo:m3}. The compact set $\mathcal{X}$ is defined by $\mathcal{X}=\{x\in\mathbb{R}:x\in[-1,1]\}$. Moreover, the cost function $J(x)$ is
\begin{equation}
    \label{eq:J_num_1}
    J(x)=-(x+0.6)^2+2.
\end{equation}
The constraint function $h(x,\delta)$ is
\begin{equation}
    h(x,\delta)=x^2+\delta-2
\end{equation}
where $\delta\sim\mathcal{N}(m_{\delta},\Sigma_{\delta}), m_{\delta}=0$, and $\Sigma_{\delta}=1$. The probability level $\alpha$ is 0.05. The optimal solution from method \textbf{Dirac-Delta} is $x^*_{\alpha}=0.595$ and the optimal objective value is $0.572$, which is plotted in Fig. \ref{fig:Result_example_one} (a). In \textbf{Dirac-Delta}, we set $\epsilon=\alpha$, $N=2000$, and $\gamma=0.01$. Besides, Fig. \ref{fig:Result_example_one} (b) and (c) show the discrete measure obtained by \textbf{Sample} and the probability density function obtained by \textbf{GMM}, respectively. For \textbf{Sample}, we choose samples $-1, -0.98, -0.96, ..., 0.96, 0.98,1$ from $\mathcal{X}$ ($S=201$) and $2000$ randomly extracted samples from $\Delta$ ($N=2000$). For \textbf{GMM}, we extracted 2000 samples from $\Delta$ randomly. Besides, we choose $L=6$. The solutions of \textbf{Sample} and \textbf{GMM} satisfy the chance constraints. For the objective function, \textbf{Sample} achieves $0.5601$ and \textbf{GMM} achieves $0.5615$, which are better than the optimal objective value achieved by \textbf{Dirac-Delta}.

A more comprehensive analysis of CPU time and sample numbers is summarized in Table 1. The CPU time increases as the sample size increases for each method. Unsurprisingly, \textbf{Sample} has a very fast computation time since it only needs to solve a linear program. In this example, since it is one dimension, the required sample number for obtaining good samples in \textbf{Sample} or approximating probability integration in \textbf{GMM} is few. It can achieve acceptable accuracy with only 50 samples. However, if the dimension of $x$ increases, the "Curse of Dimensionality" will emerge. We will show it in the second example.

\subsection{Quadrotor System Control}
\begin{figure}
\centering
\includegraphics[scale=0.375]{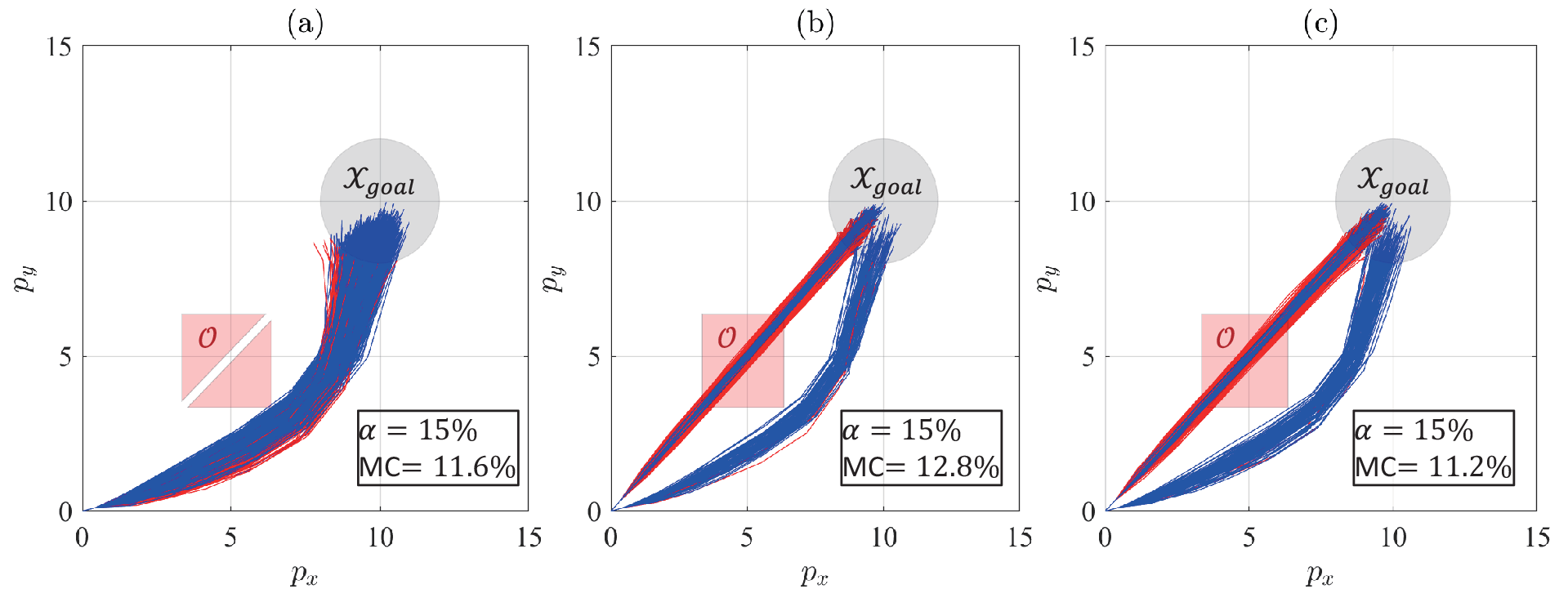}
\centering
\caption{Solutions from different methods for the tolerable failure probability threshold $\alpha=15\%$. Blue trajectories from Monte-Carlo (MC) simulations denote feasible trajectories that reach the goal set $\mathcal{X}_{goal}$ and avoid obstacles $\mathcal{O}$. Red trajectories violate constraints: (a) \textbf{Dirac-Delta} ($\text{MC}=11.6\%$ represents that the violation probability is $11.6\%$ in the MC simulations); (b) \textbf{Sample} ($\text{MC}=12.8\%$ represents that the violation probability is $12.8\%$ in the MC simulations); (c) \textbf{GMM} ($\text{MC}=11.2\%$ represents that the violation probability is $11.2\%$ in the MC simulations).}
\label{fig:control_exa}
\end{figure}

The second example considers a quadrotor system control problem in turbulent conditions. The control problem is expressed as follows:
\begin{equation} \tag{$P_{\text{QSC}}$}
\begin{split} 
&\underset{\mu\in M(\mathcal{U}^T)}{{\normalfont\text{min}}} \,\, \mathbb{E}\{\ell^x(x)+\ell^u(u)\} \\
&{\normalfont \text{s.t.}}\quad  x_{t+1}=Ax_t+B(m)u_t+d(x_t,\varphi)+\omega_t,\ u\sim M(\mathcal{U}^T), \\
& \ \ \ \ \quad t=0,1,...,T-1,\\
&\ \ \ \ \quad \mathsf{Pr}\{\left(x_t\notin\mathcal{O},\forall t=1,...,T-1\right),\ (x_T\in\mathcal{X}_{\text{goal}})\} \geq 1-\alpha,
\end{split}
\end{equation}
where $A$, $B(m)$, $d(x_t,\varphi)$ are written by
\begin{equation*}
    A=
    \begin{bmatrix}
    1 & \Delta t & 0 & 0 \\
    0 & 1 & 0 & 0 \\
    0 & 0 & 1 & \Delta \\
    0 & 0 & 0 & 1
    \end{bmatrix},\ 
    B(m)=\frac{1}{m}
    \begin{bmatrix}
    \frac{\Delta t^2}{2} & 0 \\
    \Delta t & 0 \\
    0 & \frac{\Delta t^2}{2} \\
    0 & \Delta
    \end{bmatrix},\ 
    d(x_t,\varphi)=-\varphi
    \begin{bmatrix}
    \frac{\Delta t^2|v_x|v_x}{2} \\
    \Delta t|v_x|v_x \\
    \frac{\Delta t^2|v_y|v_y}{2} \\
    \Delta t|v_y|v_y
    \end{bmatrix},
\end{equation*}
and $\Delta t$ is the sampling time, the state of the system is denoted as $x_t=[p_{x,t},v_{x,t},p_{y,t},v_{y,t}]\in\mathbb{R}^4$, the control input of the system is $u_t=\{u_{x,t},u_{y,t}\}$ within $\mathcal{U}:=\{u_t\in\mathbb{R}^2:-10\leq u_{x,t}\leq 10,-10\leq u_{y_t}\leq 10\}$, and the state and control trajectories are denoted as $x=(x_t)_{t=1}^{T}$ and $u=(u_t)_{t=1}^{T-1}$. The system starts from an initial point $x_0=[-0.5,0,-0.5,0]$. The system is expected to reach the destination set $\mathcal{X}_{\text{goal}}=\{x\in\mathbb{R}^4|\|(p_x-10,p_y-10)\|\leq 2\}$ at time $T=10$ while avoiding two polytopic obstacles $\mathcal{O}$ shown in Fig. \ref{fig:control_exa}. $\mathcal{O}$ is defined by the following constraints:
\begin{equation*}
    p_{x,t}\leq 6.35,\ p_{y,t}\geq 3.35,\ p_{x,t}-0.2-p_{y,t}\geq 0,
\end{equation*}
\begin{equation*}
    p_{x,t}\geq 3.35,\ p_{y,t}\leq 6.35,\ p_{x,t}+0.2-p_{y,t}\leq 0.
\end{equation*}
The dynamics are parametrized by uncertain parameter vector $\delta_t=[m,\varphi]^\top$, where $m>0$ represents the system's mass and $\varphi>0$ is an uncertain drag coefficient. The parameter vector $\delta$ of the system is uncorrelated random variables such that $(m-0.75)/0.5\sim \text{Beta}(2,2)$ and $(\varphi-0.4)/0.2\sim \text{Beta}(2,5)$, where $\text{Beta}(a,b)$ denotes a Beta distribution with shape parameters $(a,b)$. $\omega_t\in\mathbb{R}^4$ is the uncertain disturbance at time step $t$, which obeys multivariate normal distribution with zero means and covariance matrix
\begin{equation*}
    \Sigma=
    \begin{bmatrix}
    0.01 & 0 & 0 & 0 \\
    0 & 0.75 & 0 & 0 \\
    0 & 0 & 0.01 & 0 \\
    0 & 0 & 0 & 0.75
    \end{bmatrix}.
\end{equation*}
For the cost function, we adopt
\begin{equation*}
    \ell^x(x) = \frac{1}{T}\sum_{t=0}^{T-1}\left((p_{x,t+1}-p_{x,t})^2+(p_{y,t+1}-p_{y,t})^2\right),
\end{equation*}
\begin{equation*}
    \ell^u(u) = \frac{0.1}{T}\sum_{t=0}^{T-1}\left(u_{1,t}^2+u_{2,t}^2\right).
\end{equation*}

Results are shown in Fig. \ref{fig:control_exa} for different methods by setting $\alpha$ as $15\%$. Fig. \ref{fig:control_exa} shows 5,000 Monte-Carlo (MC) simulations of the quadrotor system using the open-loop policy computed using \textbf{Dirac-Delta} ($\epsilon=\alpha,\gamma=0.01, N=2000$), \textbf{Sample} ($S=5.1\times 10^6,N=2000$), and \textbf{GMM} ($L=6,N=2000$). When using \textbf{Dirac-Delta}, the algorithm gives a deterministic control policy that satisfies the desired success probability $1-\alpha$. When using \textbf{Sample}, or \textbf{GMM}, the algorithm gives a stochastic control policy that satisfies the desired success probability $1-\alpha$. The control inputs that generate trajectories passing through the riskier middle corridor between the obstacles are selected randomly for the stochastic control policies. The costs by using \textbf{Sample} and \textbf{GMM} are reduced by $8.2\%$ and $7.9\%$ compared to using \textbf{Dirac-Delta}. This shows that our approach can compute a better policy that solves the problem than a deterministic policy.  

\begin{figure}
\centering
\includegraphics[scale=0.425]{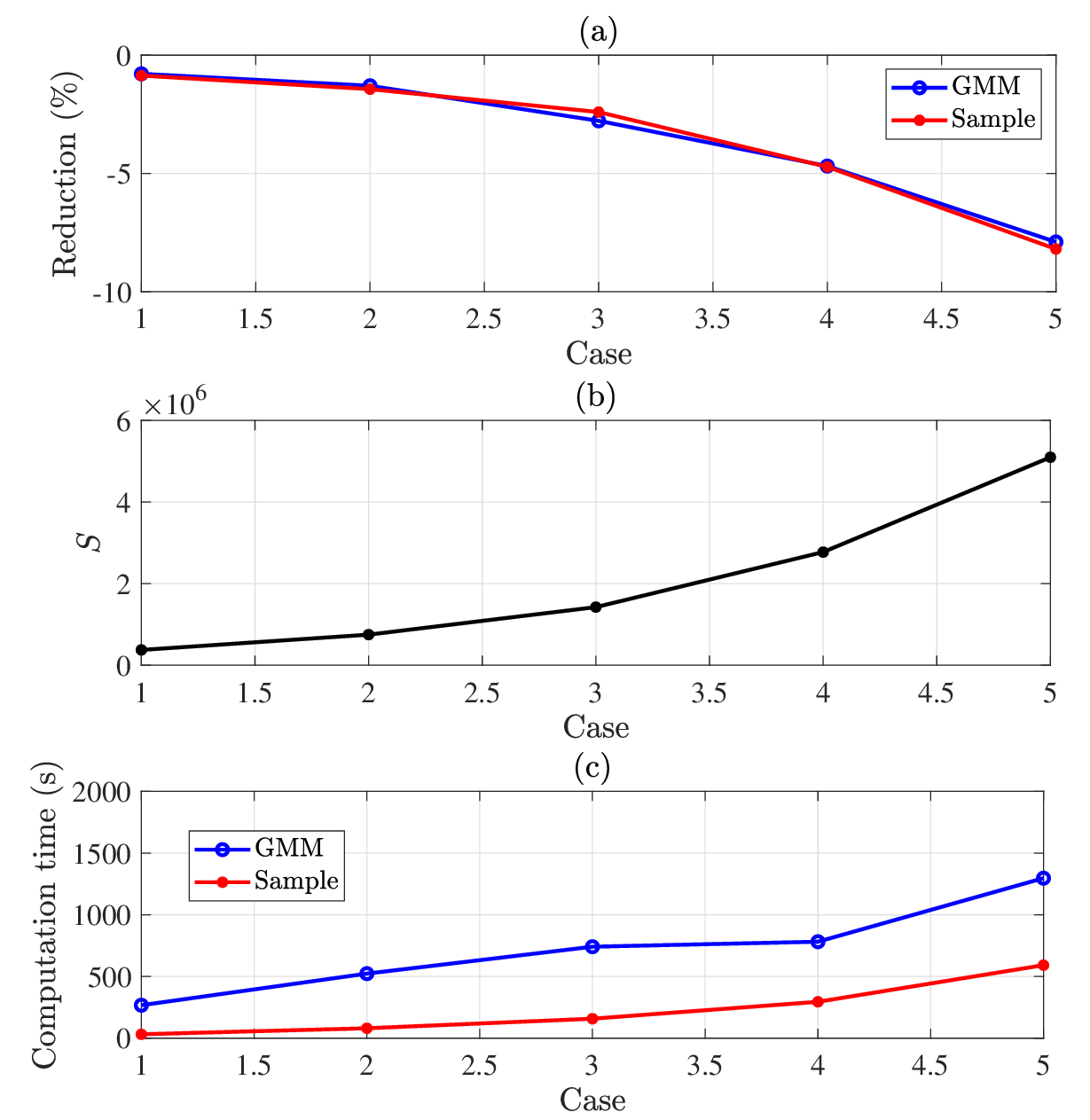}
\centering
\caption{The statistics of the control performance: (a) Reduction of cost; (b) Required samples; (c) Computation time.}
\label{fig:statistics_results}
\end{figure}
A more comprehensive comparison between the GMM-based and sample-based approximations is plotted in Figure \ref{fig:statistics_results}. Five cases are considered with different sample numbers for extracting the control input. Figure \ref{fig:statistics_results} (a) shows that the two algorithms similarly reduce the optimal objective function value. Figure \ref{fig:statistics_results} (b) shows each case's used sample number $S$ of decision variables. By comparing Figure \ref{fig:statistics_results} (a) and (b), we can see that enough samples are required to ensure the performance of the approximations. As shown in Figure \ref{fig:statistics_results} (c), the computation time increases dramatically as the sample number increases. In this comparison, for \textbf{GMM}, we choose $L=6$, and the probability integration is approximated by using the same samples of \textbf{Sample}. The computation time of \textbf{GMM} is even longer than \textbf{Sample}. One way to decrease the computation time of \textbf{GMM} is to develop fast algorithms for probability integration. We leave this for future work. In this example, the dimension of the decision variable is $20$. If the dimension increases, the required sample number will increase, and the computation time will consequently increase for \textbf{Sample} and \textbf{GMM}. We leave the issue of the "Curse of Dimensionality" for future work.  

\section{Conclusions}
\label{sec:conclusions}

In conclusion, the chance-constrained linear program in probability measure space has been addressed using sample approximation or function approximation. We establish optimization problems in finite vector space as approximate problems of chance-constrained linear programs in probability measure space. By solving the approximate problems, we can obtain the approximate solution of the chance-constrained linear program in probability measure space. Numerical examples have been implemented to validate the performance of the proposed method. Future work will be focused on the following points:
\begin{itemize}
    \item To implement sample approximation method $\tilde{P}_{\alpha}(\mathcal{C}_S,\mathcal{D}_N)$, samples of decision variable are required. As the dimension of the decision variable increases, the required sample number for a good approximation will also increase, bringing the issue of the "Curse of dimensionality." To overcome the issue of the "Curse of Dimensionality," it is important to develop efficient sampling algorithms to get "good but small samples" to ensure good approximation performance and mitigate the computation burden;
    \item For Gaussian mixture model-based approximation method $\hat{P}_{\alpha}(L,\mathcal{D}_N)$, the remaining issue is how to approximate the probability integration by fast algorithms when the problem is with complex cost function and constrained functions in high dimension space.
\end{itemize}

\begin{acknowledgements}
We thank two anonymous reviewers for taking the precious time and effort to review our manuscript and give us valuable suggestions.
\end{acknowledgements}

%References


\begin{thebibliography}{plain}

\bibitem{Ackooij} 
van Ackooij, W., Henrion, R.: Gradient formulate for nonlinear probabilistic constraints with Gaussian and Gaussian-like distributions. SIAM Journal on Optimization 24, 1864-1889 (2014)

\bibitem{Ackooij:2019} 
van Ackooij, W., Perez-Aros, P.: Generalized differentiation of probability functions acting on an infinite system of constraints. SIAM Journal on Optimization 29(3), 2179-2210 (2019)

\bibitem{Ackooij:2020} 
van Ackooij, W., Henrion, R., Perez-Aros, P.: Generalized gradients for probabilistic/robust (probust) constraints. Optimization 69(7-8), 1451-1479 (2020)

\bibitem{Berthold:2022} 
Berthold, H., Heitsch, H., Henrion, R., Schwienteck, J.: On the algorithmic solution of optimization problems subject to probabilistic/robust(probust) constraints. Mathematical Methods of Operations Research 96, 1-37 (2022)

\bibitem{Blackmore:2011} 
Blackmore, L., Ono, M., Williams, B.C.: Chance-constrained optimal path planning with obstacles. IEEE Transactions on Robotics 27(6), 1080-1094 (2011)

\bibitem{Bertsekas}  
Bertsekas, D.P., Tsitsiklis, J.N.: Introduction to Probability, Athena Scientific, Belmont, Massachusetts  (2002)

\bibitem{Billingsley}  
Billingsley, P.: Probability and Measure, John Wiley $\&$ Sons, New York  (1995)

\bibitem{Calariore:2006} 
Calariore, G., Campi, M.C.: The scenario approach to robust control design. IEEE Transactions on Automatic Control 51(5), 742-753 (2006)

\bibitem{Campi:2008} 
Campi, M.C., Garatti, S.: The exact feasibility of randomized solutions of uncertain convex programs. SIAM Journal on Optimization 19(3), 1211-1230 (2008)

\bibitem{Campi:2011} 
Campi, M.C., Garatti, S.: A sampling-and-discarding approach to chance-constrained optimization: feasibility and optimality. Journal of Optimization Theory and Applications 148(2), 257-280 (2011)

\bibitem{Campi_unconvex} 
Campi, M.C., Garatti, S.: A general scenario theory for nonconvex optimization and decision making. IEEE Transactions on Automatic Control 63(12), 4067-4078 (2015)

\bibitem{CampiBook}  
Campi, M.C., Garatti, S., Ramponi, F.A.: Introduction to the scenario approach, MOS-SIAM Series on Optimization, Philadelphia  (2019)

\bibitem{Castillo:2019} 
Castillo-Lopez, M., Ludivig, P., Sajadi-Alamdari, S.A.: A real-time approach for chance-constrained motion planning with dynamic obstacles. IEEE Robotics and Automation Letters 5(2), 3620-3625 (2020)

\bibitem{Chen:2021} 
Chen, P., Ghattas, O.: Taylor approximation for chance constrained optimization problems governed by partial differential equations with high-dimensional random parameters. SIAM/ASA Journal on Uncertainty Quantification 9(4), 1381-1410 (2021)

\bibitem{Geletu:2019} 
Geletu, A., Hoffmann, A., Kloppel, M., Li, P.: An inner-outer approximation approach to chance constrained optimization. SIAM Journal on Optimization 27(3), 1834-1857 (2017)

\bibitem{Geletu:2020} 
Geletu, A., Hoffmann, A., Schmidt, P., Li, P.: Chance constrained optimization of elliptic PDE systems with a smoothing convex approximation. ESAIM: Control, Optimisation and Calculus of Variations 26(70), (2020)

\bibitem{Grandon:2022} 
Grandon, T.G., Henrion, R., Perez-Aros, P.: Dynamic probabilistic constraints under continuous random distributions. Mathematical Programming 196, 1065-1096 (2022)

\bibitem{Gryazina} 
Gryazina, E., Polyak, B.: Random sampling: Billiard walk algorithm. European Journal of Operational Research 238, 497-504 (2014)

\bibitem{Hewing:2020} 
Hewing, L., Kabzan, J., Zeilinger, M.N.: Cautious model predictive control using Gaussian process regression. IEEE Transactions on Control Systems Technology 28(6), 2736-2743 (2020)

\bibitem{Kibzun}  
Kibzun, A.I., Kan, Y.S.,: Stochastic Programming Problems, Wiley, West Sussex, Engand  (1996)

\bibitem{Lew:2022} 
Lew, T., Sharma, A., Harrison, J., Bylard, A., Pavone, M.: Safe active dynamics learning and control: asequential exploration-exploitation framework. IEEE Transactions on Robotics 38(5), 2888-2907 (2022)

\bibitem{Luedtke:2008} 
Luedtke, J., Ahmed, S.: A sample approximation approach for optimization with probabilistic constraints. SIAM Journal on Optimization 19(2), 674-699 (2008)

\bibitem{Luedtke:2010} 
Luedtke, J., Ahmed, S., Nemhauser, G.L.: An integer programming approach for linear programs with probabilistic constraints. Mathematical Programming 122, 247-272 (2010)

\bibitem{Molchanov} 
Molchanov, I.: Theory of Random Sets. Springer, London (2005)

\bibitem{Nemirovski} 
Nemirovski, A., Shapiro, A.: Convex approximations of chance constrained programs. SIAM Journal on Optimization 17, 969-996 (2007)

\bibitem{Ono:2015} 
Ono, M., Pavone, M., Kuwata, Y., Balaram, J.: Chance-constrained dynamic programming with application to risk-aware robotic space exploration. Autonomous Robots 39(4), 555 - 571 (2015)

\bibitem{Pagnoncelli:2009} 
Pagnoncelli, B.K., Ahmed, S., Shapiro, A.: Computational study of a chance constrained portfolio selection problem. Journal of Optimization Theory and Applications 142(2), 399 - 416 (2009)

\bibitem{Pena:2020} 
Pena-Ordieres, A., Luedtke, J., Wachter, A.: Solving chance-constrained problems via a smooth sample-based nonlinear approximation. SIAM Journal on Optimization 30(3), 2221-2250 (2020)

\bibitem{Shaker:2018} 
Farshbaf-Shaker, M.H., Henrion, R., Homber, D.: Properties of chance constraints in infinite dimensions with an application to PDE constrained optimization. Set-Valued and Variational Analysis 26, 821-841 (2018)

\bibitem{Shapiro}  
Shapiro, A., Dentcheva, D., Ruszczynski, A.: Lectures on Stochastic Programming: Modeling and Theory, 2nd ed., SIAM, Philadelphia  (2014)


\bibitem{Shen_ETCI} 
Shen, X., Ouyang, T., Zhang, Y., Zhang, X.: Computing probabilistic bounds on state trajectories for uncertain systems. IEEE Transactions on Emerging Topics in Computational Intelligence 7(1), 285 - 290 (2023)

\bibitem{Shen:2021} 
Shen, X., Ouyang, T., Yang, N., Zhuang, J.: Sample-based neural approximation approach for probabilistic constrained programs. IEEE Transactions on Neural Networks and Learning Systems 34(2), 1058 - 1065 (2023)

\bibitem{Smith} 
Smith, R.L.: Efficient Monte Carlo procedures for generating points uniformly distributed over bounded regions. Operations Research 32(6), 1296 - 1308 (1984)

\bibitem{Teo:2013} 
Wu, C., Teo, K.L., Wu, S.: Min-max optimal control of linear systems with uncertainty and terminal state constraints. Automatica 49(6), 1809 - 1815 (2013)

\bibitem{Teo:2017} 
Sun, Y., Aw, G., Loxton, R., Teo, K.L.: Chance-constrained optimization for pension fund portfolios in the presense of default risk. European Journal of Operational Research 256(1), 205-214 (2017)

\bibitem{Thorpe:2022} 
Thorpe, A.J., Lew, T., Oishi, M.M.K, Pavone, M.: Data-driven chance constrained control using kernel distribution embeddings. Proceedings of Machine Learning Research 144, 1 - 13 (2022)


\bibitem{Zeevi:1997} 
Zeevi, A.J., Meir, R..: Density estimation through convex combinations of densities: approximation and estimation bounds. SIAM Journal on Optimization 10(1), 99-109 (1997)

\end{thebibliography}
\end{document}